\documentclass[11pt]{amsart}
\usepackage[margin=1.1in]{geometry}
\usepackage{amsmath,amssymb,amsthm,mathrsfs,mathtools}
\usepackage{tikz-cd}
\usepackage{enumitem}
\usepackage{booktabs,array}
\usepackage{xcolor}
\usepackage[colorlinks=true,linkcolor=blue!50!black,citecolor=green!40!black,urlcolor=blue!50!black]{hyperref}

\numberwithin{equation}{section}
\makeatletter
\providecommand{\etalchar}[1]{$^{#1}$}
\def\subsection{\@startsection{subsection}{2}{\z@}{.7\linespacing\@plus\linespacing}{.5\linespacing}{\normalfont\bfseries}}
\def\subsubsection{\@startsection{subsubsection}{3}{\z@}{.7\linespacing\@plus\linespacing}{.5\linespacing}{\normalfont\itshape}}
\makeatother

\newtheorem{theorem}{Theorem}[section]
\newtheorem{proposition}[theorem]{Proposition}
\newtheorem{lemma}[theorem]{Lemma}
\newtheorem{corollary}[theorem]{Corollary}
\newtheorem{question}[theorem]{Question}
\theoremstyle{definition}
\newtheorem{definition}[theorem]{Definition}

\theoremstyle{remark}
\newtheorem{remark}[theorem]{Remark}

\newcommand{\Ku}{\mathcal{K}u}
\newcommand{\CC}{\mathbb{C}}
\newcommand{\PP}{\mathbb{P}}
\newcommand{\cO}{\mathcal{O}}
\newcommand{\cU}{\mathcal{U}}
\newcommand{\cQ}{\mathcal{Q}}
\newcommand{\cP}{\mathcal{P}}
\newcommand{\cE}{\mathcal{E}}
\newcommand{\cF}{\mathcal{F}}
\newcommand{\cB}{\mathcal{B}}
\newcommand{\cC}{\mathcal{C}}

\newcommand{\Db}{\mathrm{D}^{\mathrm{b}}}
\DeclareMathOperator{\Hom}{Hom}

\DeclareMathOperator{\Ext}{Ext}
\DeclareMathOperator{\Ker}{Ker}

\DeclareMathOperator{\im}{Im}
\DeclareMathOperator{\Aut}{Aut}
\DeclareMathOperator{\Gr}{Gr}
\DeclareMathOperator{\HH}{HH}
\DeclareMathOperator{\NHH}{NHH}
\DeclareMathOperator{\Sym}{Sym}
\DeclareMathOperator{\res}{res}

\DeclareMathOperator{\Pic}{Pic}

\DeclareMathOperator{\ch}{ch}
\DeclareMathOperator{\Spec}{Spec}

\DeclareMathOperator{\PGL}{PGL}

\begin{document}

\title[Infinitesimal Torelli problems for special Fano threefolds]{Infinitesimal Torelli problems for special Gushel--Mukai and related Fano threefolds: Hodge-theoretic and categorical perspectives}

\author{Xun Lin}
\address{School of Science and Engineering, The Chinese University of Hong Kong, Shenzhen, Shenzhen, Guangdong 518172, China}
\email{linxun@cuhk.edu.cn}
\author{Shizhuo Zhang}
\address{School of Mathematics, Sun Yat-sen University, Guangzhou 510275, P.~R.~China}
\email{zhangshzh28@mail.sysu.edu.cn}
\author{Zheng Zhang}
\address{Institute of Mathematical Sciences, ShanghaiTech University, Shanghai 201210, China}
\email{zhangzheng@shanghaitech.edu.cn}

\date{October 16, 2026}

\begin{abstract}
A hyperelliptic curve is the classical example of a variety for which the global Torelli theorem holds while the differential of the period map has a kernel, forced by the extra involution. We study the analogous phenomenon for Fano threefolds that are double covers of rigid Fano threefolds branched along a K3 surface, mainly special Gushel--Mukai threefolds. We prove that the part of the infinitesimal period map invariant under the covering involution is injective, so that the kernel is exactly the three-dimensional space of anti-invariant deformations, one more than for ordinary Gushel--Mukai threefolds. We give a Hodge-theoretic proof, which reduces to a twisted infinitesimal Torelli theorem for the branch K3 surface, and a categorical proof, based on Kuznetsov's normal Hochschild cohomology, which also applies in genus eight and identifies the kernel with the tangent space to a distinguished point of a Bridgeland moduli space. The two descriptions of the kernel are related directly through the Pl\"ucker net and the conic of special lines of the quintic del Pezzo threefold. Finally, ordinary Verra threefolds satisfy the infinitesimal Torelli theorem, while special ones have a one-dimensional kernel; in both cases the kernel agrees with that of the categorical period map of a Clifford component, and it is the tangent space to the deformations of a canonical spinor bundle, which is rigid in the ordinary case and obstructed at second order in the special case.
\end{abstract}

\maketitle
\setcounter{tocdepth}{1}
\tableofcontents
\newpage

\section{Introduction}

\subsection{Torelli problems and the hyperelliptic phenomenon}
The Torelli problem asks whether a variety is determined by the Hodge structure on its cohomology. Its prototype is the theorem for curves: the period map $\cP_g\colon M_g\to A_g$, sending a curve of genus $g\ge2$ to its principally polarized Jacobian, is injective. The \emph{infinitesimal} Torelli problem asks whether the differential of the period map is injective. For curves this differential is the cup product map
\[
d\cP_g\colon H^1(C,T_C)\longrightarrow\Hom\bigl(H^0(C,\omega_C),H^1(C,\cO_C)\bigr),
\]
and by Noether's theorem it is injective precisely when $C$ is not hyperelliptic. On the hyperelliptic locus the differential has a $(g-2)$-dimensional kernel, although the global Torelli theorem still holds. The reason is the hyperelliptic involution: deformations anti-invariant under it are invisible to first order in the Hodge structure \cite{Usu81,Cat84}.

Fano threefolds provide a rich laboratory for Torelli problems. Their intermediate Jacobians played a central role in the work of Clemens and Griffiths on the irrationality of cubic threefolds, and in the study of prime Fano threefolds of genus six, the \emph{Gushel--Mukai threefolds} \cite{Ili94,Log12,DIM12,KP18,Deb20}. Gushel--Mukai threefolds come in two types. An \emph{ordinary} one is the intersection of the Grassmannian $\Gr(2,5)\subset\PP^9$ with a linear space of codimension two and a quadric. A \emph{special} one is a double cover of the quintic del Pezzo threefold $Y=\Gr(2,5)\cap\PP^6$ branched along a K3 surface $S\in|-K_Y|$. Special Gushel--Mukai threefolds carry an extra involution and form a closed subfamily of codimension three in the $22$-dimensional moduli space; they are the Gushel--Mukai analogue of hyperelliptic curves.\footnote{This should not be confused with the notion of a hyperelliptic Fano threefold.} For ordinary Gushel--Mukai threefolds, Debarre, Iliev and Manivel proved that the differential of the period map always has a two-dimensional kernel \cite{DIM12}; the period map is not even generically injective on the moduli space. One of our goals is to determine what happens on the special locus.

A second, more recent viewpoint replaces Hodge structures by derived categories. A prime Fano threefold $X$ of genus $g\ge6$ has a semiorthogonal decomposition
\[
\Db(X)=\langle\Ku(X),\cO_X,\cU_X^\vee\rangle,
\]
where $\cU_X$ is the restriction of a tautological bundle and $\Ku(X)$ is its \emph{Kuznetsov component} \cite{Kuz09b}. The Kuznetsov component is expected to carry the essential geometry of $X$, and its Hochschild (co)homology recovers the relevant Hodge-theoretic data. In this language there is an \emph{infinitesimal categorical period map}
\[
\eta\colon H^1(X,T_X)\longrightarrow\HH^2(\Ku(X)),
\]
sending a first-order deformation of $X$ to the induced first-order deformation of its Kuznetsov component \cite{JLLZ23}. A natural question is how the two kinds of period maps compare, and what geometric deformation problem their kernels represent.

We also study \emph{Verra threefolds}, which occur in the theory of Prym varieties of plane sextics and led to a counterexample to the tetragonal conjecture \cite{Ver04}. As for Gushel--Mukai threefolds, there are ordinary Verra threefolds, smooth divisors of bidegree $(2,2)$ in $\PP^2\times\PP^2$, and special ones, double covers of a divisor of bidegree $(1,1)$ branched along an anticanonical K3 surface. We apply both the classical Hodge-theoretic method and the categorical method to the infinitesimal Torelli problem for Verra threefolds, and we give a geometric interpretation of the kernel of the differential of the period map: for a fixed conic bundle structure $f\colon X\to\PP^2$, it is the tangent space, at the canonical spinor bundle, to the local moduli stack of spinor bundles on the conic bundle.

\subsection{Main results}
Throughout, the infinitesimal period map of a smooth Fano threefold $X$ is the map induced by cup product,
\[
d\cP\colon H^1(X,T_X)\longrightarrow\Hom\bigl(H^1(X,\Omega^2_X),H^2(X,\Omega^1_X)\bigr),
\]
and we say that $X$ satisfies the infinitesimal Torelli theorem if $d\cP$ is injective.

Let $X\to Y$ be a special Gushel--Mukai threefold, with branch surface $S\in|-K_Y|$ and covering involution $\iota$. By a lemma of Konno \cite{Kon85}, the deformation space of $X$ splits into an invariant and an anti-invariant part,
\[
H^1(X,T_X)\simeq H^1\bigl(Y,T_Y(-\log S)\bigr)\oplus H^1\bigl(Y,T_Y(-1)\bigr),
\]
of dimensions $19$ and $3$. The first summand consists of the deformations of $X$ that remain special, which, since $Y$ is rigid, are the same as the deformations of $S$ inside $Y$; the second consists of the deformations in the direction normal to the special locus.

\begin{theorem}[Theorems~\ref{thm:kernelGM} and~\ref{thm:invariantGM}]\label{thm:main1}
Let $X$ be a special Gushel--Mukai threefold. Then the invariant part of the infinitesimal period map,
\[
H^1\bigl(Y,T_Y(-\log S)\bigr)\longrightarrow\Hom\Bigl(H^1\bigl(Y,\Omega^2_Y(\log S)(-1)\bigr),H^2\bigl(Y,\Omega^1_Y(\log S)(-1)\bigr)\Bigr),
\]
is injective. Consequently $\Ker(d\cP)=H^1(Y,T_Y(-1))$, which has dimension three.
\end{theorem}

Thus every deformation of $X$ that keeps it special is detected, to first order, by the Hodge structure, and the kernel of $d\cP$ consists precisely of the three directions normal to the special locus, which the involution renders invisible. This is exactly what happens for hyperelliptic curves, where the kernel consists of the anti-invariant deformations. Note that the kernel jumps from dimension two at ordinary Gushel--Mukai threefolds to dimension three on the special locus.

For a prime Fano threefold $X$ of genus $g\ge6$, the infinitesimal categorical period map $\eta\colon H^1(X,T_X)\to\HH^2(\Ku(X))$ is the categorical analogue of the Kodaira--Spencer map: it sends a first-order deformation of $X$ to the first-order deformation of $\Ku(X)$ that it induces.

\begin{theorem}[Theorem~\ref{thm:categorical} and Corollary~\ref{cor:dimensions}]\label{thm:main2}
Let $X$ be a prime Fano threefold of genus six or eight, and let $\cQ_X$ be the tautological quotient bundle. Then
\[
\Ker(d\cP)=\Ker(\eta)\simeq\Hom(\cU_X,\cQ_X^\vee).
\]
In particular $\Ker(d\cP)$ has dimension $2$ for ordinary and $3$ for special Gushel--Mukai threefolds, and dimension $5$ in genus eight.
\end{theorem}

The theorem says two things. A first-order deformation of $X$ is invisible to the Hodge structure exactly when it leaves the Kuznetsov component unchanged, so the kernel of $d\cP$ is a categorical invariant; and this kernel has a simple description in terms of the tautological bundles. For special Gushel--Mukai threefolds we obtain a second proof of Theorem~\ref{thm:main1}. The analogous results in the other genera are recalled in Remark~\ref{rem:othergenera}.

For a general ordinary Gushel--Mukai threefold, the connected component through $X$ of the fibre of the period map is the minimal model of its Fano surface of conics, and this surface is a Bridgeland moduli space $M_\sigma(\Ku(X),[I_C])$ \cite{DIM12,JLLZ24}. Hence $\Ker(d\cP)$ is the tangent space of this moduli space at the point corresponding to $X$. The following theorem extends this picture.

\begin{theorem}[Theorem~\ref{thm:Bridgeland}]\label{thm:main3}
Let $X$ be a prime Fano threefold of genus six or eight. Then $\PP\Ker(d\cP)\simeq\PP\Ker(\eta)\simeq\PP\Hom(\cU_X,\cQ_X^\vee)$ is the exceptional locus of the birational morphism $C(X)\to M_\sigma(\Ku(X),[I_C])$ from the Hilbert scheme of conics if $g=6$, and of the birational morphism $M_X^{ss}(2,0,4)\to M_\sigma^{ss}(\Ku(X),2[I_C])$ from the moduli space of semistable rank-two sheaves if $g=8$.
\end{theorem}

In other words, the projectivized kernel can be seen geometrically: it is the family of conics, or of rank-two sheaves, that is contracted to the distinguished point of the Bridgeland moduli space. It is a line for ordinary and a plane for special Gushel--Mukai threefolds, and a $\PP^4$ in genus eight.

For a special Gushel--Mukai threefold, the two computations therefore give isomorphic spaces
\[
H^1(Y,T_Y(-1))\simeq\Ker(d\cP)=\Ker(\eta)\simeq\Hom(\cU_X,\cQ_X^\vee),
\]
the first describing deformations normal to the special locus, the second tangent vectors to a Bridgeland moduli space. It is natural to ask:

\begin{question}\label{q:direct}
Is there a direct relation between $H^1(Y,T_Y(-1))$ and $\Hom(\cU_X,\cQ_X^\vee)$?
\end{question}

\begin{theorem}[Theorem~\ref{thm:polarity}]\label{thm:main4}
Let $B\subset\wedge^2V^\vee$ be the three-dimensional net of Pl\"ucker hyperplanes cutting out $Y$ in $\Gr(2,V)$. Then $H^1(Y,T_Y(-1))\simeq B^\vee$, $\Hom(\cU_X,\cQ_X^\vee)\simeq B$, and the conic of special lines on $Y$ gives a polarity $B\simeq B^\vee$. The resulting isomorphism $H^1(Y,T_Y(-1))\simeq\Hom(\cU_X,\cQ_X^\vee)$ is $\Aut(Y)$-equivariant, canonical up to a scalar, and independent of the branch surface.
\end{theorem}

Thus the two spaces are dual realizations of the same three-dimensional space of equations of $Y$, and the equality of their dimensions is explained by the projective geometry of the quintic del Pezzo threefold alone, without reference to periods or to Hochschild cohomology.

Finally, we treat ordinary and special Verra threefolds, using both methods.

\begin{theorem}[Section~\ref{sec:Verra}]\label{thm:main5}
\leavevmode
\begin{enumerate}[label=\textup{(\arabic*)}]
\item The infinitesimal Torelli theorem holds for ordinary Verra threefolds.
\item For a special Verra threefold $X\to Y$, the invariant part of the infinitesimal period map is injective, and $\Ker(d\cP)=H^1(Y,T_Y(-1,-1))$ is one-dimensional.
\item Fix one of the two conic bundle structures $f\colon X\to\PP^2$, with Clifford component $\cC_f\subset\Db(X)$ and categorical period map $\eta_f\colon H^1(X,T_X)\to\HH^2(\cC_f)$. Then $\Ker(\eta_f)=\Ker(d\cP)$ in both cases, and this space is the tangent space to the deformations of the canonical spinor bundle of $f$.
\item The canonical spinor bundle is rigid for an ordinary Verra threefold. For a special Verra threefold its deformation space is one-dimensional with a nonzero quadratic obstruction, so its rigidified miniversal deformation ring is $\CC[[t]]/(t^2)$.
\end{enumerate}
\end{theorem}

Special Verra threefolds form a divisor in the moduli space of Verra threefolds, and, as for hyperelliptic curves, the involution forces a one-dimensional kernel in the direction normal to this divisor. By part (3), the kernels of the differentials of the classical and categorical period maps agree here as well, and the kernel is the tangent space to the deformations of a natural vector bundle. Part (4) shows that this tangent direction does not integrate: in the special case the spinor bundle is a nonreduced point of its moduli space rather than a point on a curve of distinct spinor bundles. We compare the two cases in detail in Section~\ref{ss:Verracompare}.

\subsection{Outline of the proof}
The proof of Theorem~\ref{thm:main1} begins with the observation that the anti-invariant part of $d\cP$ vanishes for parity reasons. Since $Y$ is rigid, the invariant deformations of $X$ are deformations of $S$ inside $Y$, and they form the image $H^1(S,T_S)_0\subset H^1(S,T_S)$ of the Kodaira--Spencer map. Residue maps then reduce the injectivity of the invariant part to the following twisted infinitesimal Torelli theorem for the branch surface.

\begin{proposition}[Proposition~\ref{prop:twistedTorelli}]\label{prop:main}
The cup product pairing $H^1(S,T_S)_0\otimes H^1(S,\Omega^1_S(-1))\to H^2(S,\cO_S(-1))$ is nondegenerate in the first factor.
\end{proposition}

We prove it by Flenner's method \cite{Fle86}; in the end it comes down to Noether's theorem for a canonical curve of genus six on $S$. The restriction to $H^1(S,T_S)_0$ cannot be dropped (Remark~\ref{rem:restriction}).

For Theorem~\ref{thm:main2} we use the factorization of $d\cP$ through $\eta$ from \cite{JLLZ23}, and show that the second factor is injective; for special Gushel--Mukai threefolds this follows from the equivalence of their Kuznetsov components with those of ordinary ones \cite{KP23}. Because $\HH^1(\Ku(X))$ vanishes in genera six and eight, the kernel of $\eta$ is the normal Hochschild cohomology of the exceptional pair $\langle\cO_X,\cU_X^\vee\rangle$ in the sense of Kuznetsov \cite{Kuz15}, and a spectral sequence computes it as $\Hom(\cU_X,\cQ_X^\vee)$.

Theorem~\ref{thm:main3} rests on the fact that the objects of the relevant moduli space that do not lie in $\Ku(X)$ are exactly the cokernels of morphisms $\cU_X\to\cQ_X^\vee$. The projection functor contracts all of them to a single object of $\Ku(X)$, and the tangent space at that object is again $\Hom(\cU_X,\cQ_X^\vee)$.

Theorem~\ref{thm:main4} is a computation on the Grassmannian: the Borel--Weil--Bott theorem identifies $H^1(Y,T_Y(-1))$ with $B^\vee$ and $\Hom(\cU_X,\cQ_X^\vee)$ with $B$, and the polarity of the conic of special lines is the $\Aut(Y)$-invariant form that identifies the two.

For Verra threefolds, the Hodge-theoretic statements of Theorem~\ref{thm:main5} are proved as above: by Flenner's method in the ordinary case, and by the argument for special Gushel--Mukai threefolds in the special case. On the categorical side we work with the Clifford component of a conic bundle structure. The normal Hochschild cohomology of its exceptional complement vanishes in the ordinary case. In the special case it is two-dimensional, and besides the kernel of the differential of the period map it contains a class coming from a global bivector field, which we separate using the covering involution and the vanishing of $\HH^1(\cC_f)$. The statements about spinor bundles come from the even Clifford algebra, which is $\cO\oplus\cO(-2)^{\oplus3}$ in the ordinary case and $\cO\oplus\Omega^1\oplus\Omega^2$ in the special case; in the latter, Clifford multiplication on $\Omega^1$ is the wedge product, and this produces the obstruction.

In short, the Hodge-theoretic argument computes the kernel as a cohomology group on the base, whereas the categorical argument realizes it as the tangent space to a moduli problem, of Bridgeland-stable objects or of spinor bundles. Comparing the two descriptions led us to Theorem~\ref{thm:main4}.

\subsection{Organization}
Section~\ref{sec:double} recalls the Hodge theory of double covers of rigid Fano threefolds branched along an anticanonical divisor. Section~\ref{sec:GM} proves Theorem~\ref{thm:main1} and Proposition~\ref{prop:main}. Section~\ref{sec:categorical} develops the categorical approach for prime Fano threefolds of genus six and eight: it proves Theorem~\ref{thm:main2} and Theorem~\ref{thm:main3}, and constructs the isomorphism of Theorem~\ref{thm:main4}. Section~\ref{sec:Verra} treats Verra threefolds, first the ordinary and then the special ones, and ends with a comparison of the two. We work over the field of complex numbers.

\subsection*{Acknowledgements}
We thank Arend Bayer, Olivier Debarre, Sasha Kuznetsov, Jie Liu, Ziqi Liu, Laurent Manivel, Lei Song, Dingxin Zhang and Zhiwei Zheng for useful discussions on related topics. Special thanks go to Arend Bayer for many valuable comments on our work, especially those related to Theorem~\ref{thm:main3}. We are indebted to Zhiyu Liu for allowing us to include Lemmas~\ref{lem:Liu1} and~\ref{lem:Liu2} and Proposition~\ref{prop:Liu}, which were obtained in an earlier joint project of his with the second author. We are also grateful to Enrico Fatighenti and Ariyan Javanpeykar for their interest in our work. S.Z. is supported by the Sun Yat-sen University Starting Grant No.~34000-12256019, and was supported by the Institute for Basic Science (IBS-R003-D1) and NSF grant No.~DMS-1928930 (while working at the Simons Laufer Mathematical Sciences Institute in Berkeley, California). Z.Z. is supported in part by NSFC grant 12201406. Part of this work was completed while some of the authors visited the Max Planck Institute for Mathematics, the Yau Mathematical Sciences Center at Tsinghua University, the Institute of Mathematical Sciences at ShanghaiTech University and Shandong University; we are grateful for their hospitality and support.

\section{Double covers of rigid Fano threefolds}\label{sec:double}

\subsection{Logarithmic sheaves and residues}
Let $Y$ be a smooth projective variety of dimension $n$ and $S\subset Y$ a smooth prime divisor. In local coordinates $z_1,\dots,z_n$ with $S=\{z_1=0\}$, the sheaf $T_Y(-\log S)$ is generated by $z_1\partial_{z_1},\partial_{z_2},\dots,\partial_{z_n}$, and $\Omega^m_Y(\log S)$ is generated by the forms $dz_J$ and $\frac{dz_1}{z_1}\wedge dz_I$, where $dz_I$ and $dz_J$ do not contain $dz_1$. We use the exact sequences
\begin{align}
&0\longrightarrow T_S\longrightarrow T_Y|_S\longrightarrow N_{S/Y}\longrightarrow0,\label{eq:TS}\\
&0\longrightarrow T_Y(-\log S)\longrightarrow T_Y\longrightarrow N_{S/Y}\longrightarrow0,\label{eq:logT}\\
&0\longrightarrow T_Y\otimes\cO_Y(-S)\longrightarrow T_Y(-\log S)\longrightarrow T_S\longrightarrow0,\label{eq:logTS}\\
&0\longrightarrow\Omega^m_Y\longrightarrow\Omega^m_Y(\log S)\xrightarrow{\ \res\ }\Omega^{m-1}_S\longrightarrow0.\label{eq:res}
\end{align}
In local coordinates, the map in \eqref{eq:logT} sends $\sum_if_i\partial_{z_i}$ to $f_1|_{z_1=0}\,\partial_{z_1}$, the map in \eqref{eq:logTS} sends $f_1z_1\partial_{z_1}+\sum_{i\ge2}f_i\partial_{z_i}$ to $\sum_{i\ge2}f_i|_{z_1=0}\,\partial_{z_i}$, and the residue sends $f\,\frac{dz_1}{z_1}\wedge dz_2\wedge\dots+\alpha$ to $f|_{z_1=0}\,dz_2\wedge\cdots$.

\begin{lemma}\label{lem:sign}
The diagram of multiplication maps
\[
\begin{tikzcd}[column sep=large]
T_S\otimes\Omega^{m-1}_S\arrow[r,"m_1"] & \Omega^{m-2}_S\\
T_Y(-\log S)\otimes\Omega^m_Y(\log S)\arrow[u,"\res_1"]\arrow[r,"m_2"] & \Omega^{m-1}_Y(\log S)\arrow[u,"\res_2"']
\end{tikzcd}
\]
commutes up to sign: $\res_2\circ m_2=-m_1\circ\res_1$.
\end{lemma}

\begin{proof}
The statement is local. Write $\alpha=a_1z_1\partial_{z_1}+\sum_{i\ge2}a_i\partial_{z_i}$ and $\beta=\frac{dz_1}{z_1}\wedge K_1+K_2$, with $K_1,K_2$ free of $dz_1$. Then
\[
m_2(\alpha,\beta)=a_1K_1-\frac{dz_1}{z_1}\wedge\Bigl(\sum_{i\ge2}a_i\partial_{z_i}\Bigr)\cdot K_1+\Bigl(\sum_{i\ge2}a_i\partial_{z_i}\Bigr)\cdot K_2,
\]
so $\res_2m_2(\alpha,\beta)=-(\sum_{i\ge2}a_i\partial_{z_i})\cdot K_1|_{z_1=0}$, while $m_1\res_1(\alpha,\beta)=(\sum_{i\ge2}a_i\partial_{z_i})\cdot K_1|_{z_1=0}$.
\end{proof}

\subsection{Invariant and anti-invariant parts of the period map}
Let $\pi\colon X\to Y$ be a double cover of smooth projective threefolds branched along a smooth prime divisor $S$, with $L^{\otimes2}\simeq\cO_Y(S)$. The infinitesimal period map
\[
d\cP\colon H^1(X,T_X)\longrightarrow\Hom\bigl(H^1(X,\Omega^2_X),H^2(X,\Omega^1_X)\bigr)
\]
is induced by the contraction $T_X\otimes\Omega^2_X\to\Omega^1_X$, and the infinitesimal Torelli theorem holds for $X$ if it is injective.

\begin{proposition}[{\cite[Lemma~1.1]{Kon85}}]\label{prop:Konno}
The covering involution decomposes the relevant groups into invariant and anti-invariant summands:
\begin{align*}
H^1(X,T_X)&\simeq H^1\bigl(Y,T_Y(-\log S)\bigr)\oplus H^1\bigl(Y,T_Y\otimes L^{-1}\bigr),\\
H^1(X,\Omega^2_X)&\simeq H^1(Y,\Omega^2_Y)\oplus H^1\bigl(Y,\Omega^2_Y(\log S)\otimes L^{-1}\bigr),\\
H^2(X,\Omega^1_X)&\simeq H^2(Y,\Omega^1_Y)\oplus H^2\bigl(Y,\Omega^1_Y(\log S)\otimes L^{-1}\bigr).
\end{align*}
Since cup product is compatible with these decompositions \cite[p.~609]{Kon85}, $d\cP$ is the sum of an invariant part
\begin{multline*}
H^1\bigl(Y,T_Y(-\log S)\bigr)\longrightarrow\Hom\bigl(H^1(Y,\Omega^2_Y),H^2(Y,\Omega^1_Y)\bigr)\\
\oplus\Hom\bigl(H^1(Y,\Omega^2_Y(\log S)\otimes L^{-1}),H^2(Y,\Omega^1_Y(\log S)\otimes L^{-1})\bigr)
\end{multline*}
and an anti-invariant part, which exchanges the two kinds of summands.
\end{proposition}

\begin{remark}\label{rem:parity}
If $H^1(Y,\Omega^2_Y)=0$, and hence $H^2(Y,\Omega^1_Y)=0$ by Hodge symmetry, the anti-invariant part vanishes, and the invariant part is the map
\begin{equation}\label{eq:invpart}
H^1\bigl(Y,T_Y(-\log S)\bigr)\longrightarrow\Hom\Bigl(H^1\bigl(Y,\Omega^2_Y(\log S)\otimes L^{-1}\bigr),H^2\bigl(Y,\Omega^1_Y(\log S)\otimes L^{-1}\bigr)\Bigr).
\end{equation}
In particular, the anti-invariant deformations $H^1(Y,T_Y\otimes L^{-1})$ lie in $\Ker(d\cP)$.
\end{remark}

From now on we assume that $Y$ is a \emph{rigid} Fano threefold, $H^1(Y,T_Y)=0$, with $H^1(Y,\Omega^2_Y)=0$, and that $S\in|-K_Y|$ is a smooth anticanonical K3 surface. Since $T_Y\otimes\omega_Y\simeq\Omega^2_Y$, we have
\begin{equation}\label{eq:vanishTS}
H^1\bigl(Y,T_Y\otimes\cO_Y(-S)\bigr)\simeq H^1(Y,\Omega^2_Y)=0.
\end{equation}
By rigidity, \eqref{eq:logT} gives a surjection $H^0(S,N_{S/Y})\twoheadrightarrow H^1(Y,T_Y(-\log S))$, and by \eqref{eq:vanishTS}, \eqref{eq:logTS} gives an injection $H^1(Y,T_Y(-\log S))\hookrightarrow H^1(S,T_S)$.

\begin{proposition}\label{prop:KS}
The composition $H^0(S,N_{S/Y})\twoheadrightarrow H^1(Y,T_Y(-\log S))\hookrightarrow H^1(S,T_S)$ is the Kodaira--Spencer map of the family of deformations of $S$ in $Y$, that is, the connecting map of \eqref{eq:TS}. Consequently, $H^1(Y,T_Y(-\log S))$ is identified with the image
\[
H^1(S,T_S)_0=\mathrm{KS}\bigl(H^0(S,N_{S/Y})\bigr)\subset H^1(S,T_S).
\]
\end{proposition}

\begin{proof}
Both maps are obtained by lifting a section of the normal bundle to local vector fields on $Y$ and taking differences of the lifts on overlaps, restricted to $S$. The second statement follows from the surjectivity of the first map.
\end{proof}

Geometrically, $H^1(Y,T_Y(-\log S))$ is the space of first-order deformations of the pair $(Y,S)$, which, since $Y$ is rigid, are deformations of $S$ inside $Y$. Infinitesimal variations of mixed Hodge structure for pairs consisting of a Fano threefold and an anticanonical divisor are studied in \cite{AGG24}; it would be interesting to compare them with the problems considered here.

\section{Special Gushel--Mukai threefolds}\label{sec:GM}

\subsection{Set-up and the anti-invariant part}
Let $Y=\Gr(2,5)\cap\PP^6\subset\PP^9$ be a smooth codimension-three linear section of the Grassmannian, the quintic del Pezzo threefold. It is a rigid Fano threefold with $H^1(Y,\Omega^2_Y)=0$ and $\Pic(Y)=\mathbb{Z}H$, and we write $\cO_Y(m)=\cO_Y(mH)$, so that $\omega_Y\simeq\cO_Y(-2)$. Let $S=Y\cap Q$ for a quadric $Q$, a smooth K3 surface of degree ten in $|-K_Y|$. The double cover $X\to Y$ branched along $S$, with $L=\cO_Y(1)$, is a \emph{special Gushel--Mukai threefold}. It satisfies $h^1(X,T_X)=22$, and Section~\ref{sec:double} applies.

\begin{lemma}\label{lem:antiinv}
The anti-invariant summand $H^1(Y,T_Y(-1))$ of $H^1(X,T_X)$ is contained in $\Ker(d\cP)$, and it has dimension three. The invariant summand $H^1(Y,T_Y(-\log S))$ has dimension $19$, and
\begin{equation}\label{eq:19}
0\longrightarrow H^1\bigl(Y,T_Y(-\log S)\bigr)\longrightarrow H^1(S,T_S)\longrightarrow\CC\longrightarrow0.
\end{equation}
\end{lemma}

\begin{proof}
The first assertion is Remark~\ref{rem:parity}. Explicitly, if $\beta\in H^1(Y,T_Y(-1))$ and $\alpha\in H^1(Y,\Omega^2_Y(\log S)(-1))$, the involution $\iota$ satisfies $\iota(\beta\cdot\alpha)=\iota(\beta)\cdot\iota(\alpha)$; the left side is $-\beta\cdot\alpha$ and the right side is $\beta\cdot\alpha$, so the product vanishes.

From \eqref{eq:logT} we get $H^1(S,N_{S/Y})\to H^2(Y,T_Y(-\log S))\to H^2(Y,T_Y)$. Here $H^1(S,N_{S/Y})=H^1(S,\cO_S(2))=0$ and $H^2(Y,T_Y)=H^2(Y,\Omega^2_Y\otimes\omega_Y^{-1})=0$ by Kodaira--Akizuki--Nakano vanishing, so $H^2(Y,T_Y(-\log S))=0$. Since $T_Y(-S)\simeq\Omega^2_Y$, we have $H^1(Y,T_Y(-S))=0$ and $H^2(Y,T_Y(-S))=H^{2,2}(Y)=\CC$, and \eqref{eq:logTS} gives \eqref{eq:19}. Hence the invariant summand has dimension $20-1=19$, and Proposition~\ref{prop:Konno} gives the dimension $22-19=3$ of the anti-invariant summand.
\end{proof}

\begin{theorem}\label{thm:kernelGM}
Let $X$ be a special Gushel--Mukai threefold. Then $\Ker(d\cP)=H^1(Y,T_Y(-1))$, a three-dimensional space.
\end{theorem}

By Proposition~\ref{prop:Konno} and Lemma~\ref{lem:antiinv}, Theorem~\ref{thm:kernelGM} follows from the next theorem, whose proof occupies the rest of this section.

\begin{theorem}\label{thm:invariantGM}
The invariant infinitesimal period map
\[
H^1\bigl(Y,T_Y(-\log S)\bigr)\longrightarrow\Hom\Bigl(H^1\bigl(Y,\Omega^2_Y(\log S)(-1)\bigr),H^2\bigl(Y,\Omega^1_Y(\log S)(-1)\bigr)\Bigr)
\]
is injective.
\end{theorem}

\subsection{Reduction to the branch K3 surface}
Twisting \eqref{eq:res} for $m=1,2$ by $\cO_Y(-1)$ gives
\begin{align}
&0\longrightarrow\Omega^1_Y(-1)\longrightarrow\Omega^1_Y(\log S)(-1)\longrightarrow\cO_S(-1)\longrightarrow0,\label{eq:res1}\\
&0\longrightarrow\Omega^2_Y(-1)\longrightarrow\Omega^2_Y(\log S)(-1)\longrightarrow\Omega^1_S(-1)\longrightarrow0.\label{eq:res2}
\end{align}

\begin{proposition}\label{prop:diagramS}
The diagram
\[
\begin{tikzcd}[column sep=large]
H^1(S,T_S)\otimes H^1\bigl(S,\Omega^1_S(-1)\bigr)\arrow[r,"m_1"] & H^2\bigl(S,\cO_S(-1)\bigr)\\
H^1\bigl(Y,T_Y(-\log S)\bigr)\otimes H^1\bigl(Y,\Omega^2_Y(\log S)(-1)\bigr)\arrow[u,"j\otimes\res_1"]\arrow[r,"m_2"] & H^2\bigl(Y,\Omega^1_Y(\log S)(-1)\bigr)\arrow[u,"\res_2"']
\end{tikzcd}
\]
commutes up to sign, and the map $j\colon H^1(Y,T_Y(-\log S))\to H^1(S,T_S)$ is injective with image $H^1(S,T_S)_0$.
\end{proposition}

\begin{proof}
Twist Lemma~\ref{lem:sign} for $m=2$ by $\cO_Y(-1)$ and take cohomology. The statement about $j$ is Proposition~\ref{prop:KS}.
\end{proof}

\begin{lemma}\label{lem:res1surj}
$H^2(Y,\Omega^2_Y(-1))=0$; hence $\res_1\colon H^1(Y,\Omega^2_Y(\log S)(-1))\to H^1(S,\Omega^1_S(-1))$ is surjective.
\end{lemma}

\begin{proof}
Since $\Omega^2_Y\simeq T_Y\otimes\omega_Y$ and $\omega_Y\simeq\cO_Y(-2)$, we have $\Omega^2_Y(-1)\simeq T_Y(-3)$. Twisting the normal sequence $0\to T_Y\to T_{\Gr}|_Y\to\cO_Y(1)^{\oplus3}\to0$ by $\cO_Y(-3)$ and using $H^1(Y,\cO_Y(-2))=0$, it suffices to show $H^2(Y,T_{\Gr}(-3)|_Y)=0$, where $\Gr=\Gr(2,5)$. Twisting the Koszul resolution
\[
0\to\cO_{\Gr}(-3)\to\cO_{\Gr}(-2)^{\oplus3}\to\cO_{\Gr}(-1)^{\oplus3}\to\cO_{\Gr}\to\cO_Y\to0
\]
by $T_{\Gr}(-3)$ and splitting the result into short exact sequences, the vanishing follows from
\begin{equation}\label{eq:Snow}
H^5\bigl(T_{\Gr}(-6)\bigr)=H^4\bigl(T_{\Gr}(-5)\bigr)=H^3\bigl(T_{\Gr}(-4)\bigr)=H^2\bigl(T_{\Gr}(-3)\bigr)=0.
\end{equation}
Since $\omega_{\Gr}\simeq\cO_{\Gr}(-5)$, we have $T_{\Gr}(-j)\simeq\Omega^5_{\Gr}(5-j)$. By Snow's vanishing theorem \cite[p.~171]{Sno86}, $H^p(\Gr,\Omega^q_{\Gr}(k))=0$ for $1\le k\le4$ when $q>4$, or when $2p\ge q>0$; this gives the third and fourth vanishings. The second holds since $H^4(\Gr,\Omega^5_{\Gr})=H^{5,4}(\Gr)=0$, as $H^9(\Gr,\CC)=0$. For the first, Serre duality gives $H^5(\Gr,\Omega^5_{\Gr}(-1))\simeq H^1(\Gr,\Omega^1_{\Gr}(1))^\vee=0$.
\end{proof}

\begin{lemma}\label{lem:reduction}
If the pairing $H^1(S,T_S)_0\otimes H^1(S,\Omega^1_S(-1))\to H^2(S,\cO_S(-1))$ is nondegenerate in the first factor, then Theorem~\ref{thm:invariantGM} holds.
\end{lemma}

\begin{proof}
Let $0\ne\alpha\in H^1(Y,T_Y(-\log S))$, viewed in $H^1(S,T_S)_0$. By hypothesis there is $\beta\in H^1(S,\Omega^1_S(-1))$ with $m_1(\alpha,\beta)\ne0$, and by Lemma~\ref{lem:res1surj} we can write $\beta=\res_1(\beta')$. Proposition~\ref{prop:diagramS} gives $\res_2(m_2(\alpha,\beta'))=\pm m_1(\alpha,\beta)\ne0$, so $m_2(\alpha,\beta')\ne0$.
\end{proof}

\subsection{A twisted infinitesimal Torelli theorem for the branch K3 surface}

\begin{proposition}\label{prop:twistedTorelli}
The pairing
\[
H^1(S,T_S)_0\otimes H^1\bigl(S,\Omega^1_S(-1)\bigr)\longrightarrow H^2\bigl(S,\cO_S(-1)\bigr)
\]
is nondegenerate in the first factor.
\end{proposition}

The proof adapts Flenner's method \cite[Theorem~1.1]{Fle86}, as in \cite[Theorem~5.1]{DIM12}. Since $N_{S/Y}\simeq\cO_S(2)$, the conormal sequence of $S\subset Y$ and \cite[Theorem~2.1]{Fle86} give
\begin{align}
&0\longrightarrow\cO_S(-2)\longrightarrow\Omega^1_Y|_S\longrightarrow\Omega^1_S\longrightarrow0,\label{eq:conormal}\\
&0\longrightarrow\Omega^1_S(-1)\longrightarrow\Omega^2_Y(1)|_S\longrightarrow\Omega^2_S(1)\longrightarrow0.\label{eq:Flenner}
\end{align}
Fix a nowhere vanishing $\omega\in H^0(S,\omega_S)$. Contraction with $\omega$ gives an isomorphism $\lrcorner\,\omega\colon H^1(S,T_S)\simeq H^1(S,\Omega^1_S)$, and we put
\[
I\colon H^1(S,T_S)_0\hookrightarrow H^1(S,T_S)\xrightarrow{\ \lrcorner\omega\ }H^1(S,\Omega^1_S)\longrightarrow H^2(S,\cO_S(-2)),
\]
the last map being the connecting map of \eqref{eq:conormal}. The surface $S=\Gr(2,5)\cap\PP^6\cap Q$ is a Gushel--Mukai surface, and the following is an application of \cite[Lemma~2.10]{Fle86}.

\begin{lemma}\label{lem:Flennerdiagram}
In the diagram
\[
\begin{tikzcd}[row sep=large]
H^1(S,T_S)_0\arrow[d,"I"'] & \otimes & H^1\bigl(S,\Omega^1_S(-1)\bigr)\arrow[r] & H^2\bigl(S,\cO_S(-1)\bigr)\arrow[d,"\wedge\omega","\simeq"']\\
H^2\bigl(S,\cO_S(-2)\bigr) & \otimes & H^0\bigl(S,\Omega^2_S(1)\bigr)\arrow[u,"\delta"']\arrow[r] & H^2\bigl(S,\Omega^2_S(-1)\bigr),
\end{tikzcd}
\]
where $\delta$ is the connecting map of \eqref{eq:Flenner} and the horizontal arrows are cup products, we have $m(\tau,\delta\gamma)\wedge\omega=\pm I(\tau)\cdot\gamma$ for all $\tau$ and $\gamma$.
\end{lemma}

Thus Proposition~\ref{prop:twistedTorelli} follows from the injectivity of $I$ and the nondegeneracy of the bottom row in the first factor.

\begin{lemma}\label{lem:H1restr}
The restriction $H^1(Y,\Omega^1_Y)\to H^1(S,\Omega^1_Y|_S)$ is an isomorphism; in particular $h^1(S,\Omega^1_Y|_S)=1$.
\end{lemma}

\begin{proof}
The sequence $0\to\Omega^1_Y(-S)\to\Omega^1_Y\to\Omega^1_Y|_S\to0$ reduces this to $H^1(Y,\Omega^1_Y(-S))=H^2(Y,\Omega^1_Y(-S))=0$. Since $\cO_Y(-S)\simeq\omega_Y$, Serre duality gives $H^1(Y,\Omega^1_Y\otimes\omega_Y)\simeq H^2(Y,T_Y)^\vee=0$ and $H^2(Y,\Omega^1_Y\otimes\omega_Y)\simeq H^1(Y,T_Y)^\vee=0$.
\end{proof}

\begin{lemma}\label{lem:Iinj}
The map $I$ is injective.
\end{lemma}

\begin{proof}
From \eqref{eq:TS} and \eqref{eq:conormal} we obtain exact sequences
\begin{gather*}
H^0(S,\cO_S(2))\xrightarrow{\ \mathrm{KS}\ }H^1(S,T_S)\xrightarrow{\ J_1\ }H^1(S,T_Y|_S)\longrightarrow H^1(S,\cO_S(2))=0,\\
0=H^1(S,\cO_S(-2))\longrightarrow H^1(S,\Omega^1_Y|_S)\xrightarrow{\ J_2\ }H^1(S,\Omega^1_S)\longrightarrow H^2(S,\cO_S(-2)),
\end{gather*}
which are dual to each other under Serre duality. The kernel of $H^1(S,\Omega^1_S)\to H^2(S,\cO_S(-2))$ is $J_2(H^1(S,\Omega^1_Y|_S))$, so we must show that $\tau\lrcorner\omega\notin\im J_2$ for $0\ne\tau\in H^1(S,T_S)_0=\im\mathrm{KS}$. Suppose $\tau\lrcorner\omega=J_2(a)$. Since $J_1(\tau)=0$, for every $b\in H^1(S,\Omega^1_Y|_S)$,
\[
\int_SJ_2(a)\wedge J_2(b)=\int_S(\tau\lrcorner\omega)\wedge J_2(b)=-\langle\tau,J_2(b)\rangle=-\langle J_1(\tau),b\rangle=0,
\]
where $(\tau\lrcorner\omega)\wedge\gamma=-(\tau\lrcorner\gamma)\wedge\omega$ is checked in local coordinates. By Lemma~\ref{lem:H1restr}, $J_2(a)=k_1H_S$ and $J_2(b)=k_2H_S$, where $H_S$ is the restricted polarization, and $\int_SH_S^2=10\ne0$; hence $J_2(a)=0$ and $\tau=0$.
\end{proof}

\begin{remark}\label{rem:restriction}
The restriction to $H^1(S,T_S)_0$ cannot be removed: the preimage of $J_2(H^1(S,\Omega^1_Y|_S))$ under $\lrcorner\,\omega$ is a line in $H^1(S,T_S)$ contained in the kernel of the extension of $I$ to all of $H^1(S,T_S)$.
\end{remark}

\begin{lemma}\label{lem:Noether}
The pairing $H^2(S,\cO_S(-2))\otimes H^0(S,\Omega^2_S(1))\to H^2(S,\Omega^2_S(-1))$ is nondegenerate in the first factor.
\end{lemma}

\begin{proof}
Since $\Omega^2_S\simeq\cO_S$, Serre duality reduces the claim to the surjectivity of $H^0(S,\cO_S(1))\otimes H^0(S,\cO_S(1))\to H^0(S,\cO_S(2))$. A general hyperplane section $C=S\cap H$ is a smooth curve of genus six with $\omega_C\simeq\cO_C(C)\simeq\cO_C(1)$, which is very ample, so $C$ is not hyperelliptic; see also \cite[Proposition~2.12]{DK18}. By \cite[Observation~1.3]{GP00} it suffices to prove surjectivity of $H^0(C,\omega_C)\otimes H^0(C,\omega_C)\to H^0(C,\omega_C^2)$, which is Noether's theorem for non-hyperelliptic curves.
\end{proof}

\begin{proof}[Proof of Proposition~\ref{prop:twistedTorelli} and Theorem~\ref{thm:invariantGM}]
Combine Lemmas~\ref{lem:Flennerdiagram}, \ref{lem:Iinj} and~\ref{lem:Noether} to get Proposition~\ref{prop:twistedTorelli}, and then apply Lemma~\ref{lem:reduction}.
\end{proof}

\section{The categorical approach for prime Fano threefolds of genus six and eight}\label{sec:categorical}

\subsection{The categorical period map}
Let $X$ be a prime Fano threefold of genus $g$, so that $-K_X$ generates $\Pic(X)$, $-K_X^3=2g-2$, and $2\le g\le12$, $g\ne11$. Infinitesimal Torelli theorems are known for $g=2,4,5$ and for non-hyperelliptic $X$ of genus three \cite[Theorem~1.4]{JLLZ23}. In this section we treat genera six and eight; the other genera are discussed in Remark~\ref{rem:othergenera}. For $g\ge6$ there is a decomposition
\[
\Db(X)=\langle\Ku(X),\cO_X,\cU_X^\vee\rangle,
\]
where $\cU_X$ is the pullback of the tautological subbundle of the relevant Grassmannian. By \cite[Theorem~1.1]{JLLZ23} there is a commutative diagram
\begin{equation}\label{eq:factor}
\begin{tikzcd}[column sep=large]
 & \HH^2(\Ku(X))\arrow[d,"\Gamma_X"]\\
H^1(X,T_X)\arrow[ru,"\eta"]\arrow[r,"d\cP"'] & \Hom\bigl(H^1(X,\Omega^2_X),H^2(X,\Omega^1_X)\bigr).
\end{tikzcd}
\end{equation}
By the Hochschild--Kostant--Rosenberg isomorphism, $\HH_{-1}(\Ku(X))\simeq H^1(X,\Omega^2_X)$ and $\HH_1(\Ku(X))\simeq H^2(X,\Omega^1_X)$, and $\Gamma_X$ is induced by the action
\[
\gamma_X\colon\HH^2(\Ku(X))\longrightarrow\Hom\bigl(\HH_{-1}(\Ku(X)),\HH_1(\Ku(X))\bigr)
\]
of Hochschild cohomology on Hochschild homology \cite[p.~11]{JLLZ23}.

\begin{proposition}\label{prop:gammainj}
For a prime Fano threefold of genus six or eight, the map $\gamma_X$ is injective. Consequently $\Ker(d\cP)=\Ker(\eta)$.
\end{proposition}

\begin{proof}
For ordinary Gushel--Mukai threefolds this is \cite[Theorem~4.6]{JLLZ23}. If $X$ is a special Gushel--Mukai threefold, there is an ordinary Gushel--Mukai threefold $X'$ with $\Ku(X')\simeq\Ku(X)$ \cite[Theorem~1.6]{KP23}, and the equivalence is of Fourier--Mukai type \cite[Theorem~1.3]{LPZ23}. By \cite[Theorem~4.8]{JLLZ24}, $\gamma_X$ is injective if and only if $\gamma_{X'}$ is, which proves the claim. In genus eight the injectivity is \cite[Theorem~4.10]{JLLZ23}; alternatively, $\Ku(X)\simeq\Ku(Y)$ for a cubic threefold $Y$ \cite[Theorem~3.17]{Kuz04}, and $\gamma_Y$ is injective by \cite[Theorem~4.4(2)]{JLLZ23}. The last assertion follows from \eqref{eq:factor}.
\end{proof}

\subsection{The kernels of the categorical period maps}
We use Kuznetsov's normal Hochschild cohomology \cite[Section~3]{Kuz15}. For an admissible subcategory $\cB\subset\Db(X)$ generated by an exceptional collection with orthogonal $\Ku(X)$, it fits into an exact sequence
\begin{equation}\label{eq:NHH}
\cdots\longrightarrow\HH^1(X)\longrightarrow\HH^1(\Ku(X))\longrightarrow\NHH^2(\cB,X)\longrightarrow\HH^2(X)\longrightarrow\HH^2(\Ku(X))\longrightarrow\cdots
\end{equation}
\cite[Theorem~3.3]{Kuz15}, and it is computed by a spectral sequence \cite[Proposition~3.7]{Kuz15}. For an exceptional collection $E_1,\dots,E_s$ and polarization $H=-K_X$, its first page is
\begin{multline}\label{eq:E1page}
E_1^{-p,q}=\bigoplus_{1\le a_0<\dots<a_p\le s}\ \bigoplus_{k_0+\dots+k_p=q}\Ext^{k_0}(E_{a_0},E_{a_1})\otimes\cdots\\
\cdots\otimes\Ext^{k_{p-1}}(E_{a_{p-1}},E_{a_p})\otimes\Ext^{k_p-3}(E_{a_p},E_{a_0}(H)),
\end{multline}
converging to $\NHH^{q-p}(\cB,X)$.

\begin{theorem}\label{thm:categorical}
Let $X$ be a prime Fano threefold of genus six or eight, and let $\cQ_X$ be the tautological quotient bundle. Then $\Ker(\eta)\simeq\Hom(\cU_X,\cQ_X^\vee)$.
\end{theorem}

\begin{proof}
Let $\cB_X=\langle\cO_X,\cU_X^\vee\rangle$. In genera six and eight, $\HH^2(X)\simeq H^1(X,T_X)$ \cite{BFT23}, and $\HH^1(\Ku(X))=0$, by \cite[Proposition~2.12]{KP18} for Gushel--Mukai threefolds and by \cite[Theorem~8.9]{Kuz09b} in genus eight. So \eqref{eq:NHH} becomes
\[
0\longrightarrow\NHH^2(\cB_X,X)\xrightarrow{\ i\ }H^1(X,T_X)\xrightarrow{\ \eta\ }\HH^2(\Ku(X)),
\]
and $\Ker(\eta)=\im(i)\simeq\NHH^2(\cB_X,X)$. In \eqref{eq:E1page} with $E_1=\cO_X$ and $E_2=\cU_X^\vee$, we have $p\in\{0,1\}$ and $q\ge3$, and the terms of total degree two have $(p,q)=(1,3)$. Moreover $E_1^{-2,3}=0$,
\begin{align*}
E_1^{0,3}&=\Hom(\cO_X,\cO_X(H))\oplus\Hom(\cO_X,\cU_X^\vee\otimes\cU_X^\vee),\\
E_1^{-1,3}&=\Hom(\cO_X,\cU_X^\vee)\otimes\Hom(\cU_X^\vee,\cO_X(H)),
\end{align*}
where we used $\Hom(\cU_X^\vee,\cU_X^\vee(H))\simeq\Hom(\cO_X,\cU_X\otimes\cU_X^\vee(H))\simeq\Hom(\cO_X,\cU_X^\vee\otimes\cU_X^\vee)$, as $\det\cU_X^\vee=\cO_X(H)$. The differential $d_1$ is the sum of the composition $\eta_1\colon s_1\otimes s_2\mapsto s_2\circ s_1$ and a multiplication map $\eta_2$ of sections of $\cU_X^\vee$. The map $\eta_1$ factors through $\eta_2$, via
\[
\cO_X\xrightarrow{\ s_1\ }\cU_X^\vee\xrightarrow{\ \alpha(s_2)\ }\cU_X^\vee\otimes\cU_X^\vee\simeq\cU_X\otimes\cU_X^\vee\otimes\cO_X(H)\longrightarrow\cO_X(H),
\]
where $\alpha$ is the natural identification $\Hom(\cU_X^\vee,\cO_X(H))\simeq\Hom(\cO_X,\cU_X^\vee)$ and the last map is the trace pairing; hence $\Ker(d_1)=\Ker(\eta_2)$. Applying $\Hom(\cU_X,-)$ to the tautological sequence $0\to\cQ_X^\vee\to\Hom(\cO_X,\cU_X^\vee)\otimes\cO_X\to\cU_X^\vee\to0$ gives
\[
0\longrightarrow\Hom(\cU_X,\cQ_X^\vee)\longrightarrow\Hom(\cO_X,\cU_X^\vee)^{\otimes2}\xrightarrow{\ \eta_2\ }\Hom(\cO_X,\cU_X^\vee\otimes\cU_X^\vee)\longrightarrow\Ext^1(\cU_X,\cQ_X^\vee)\longrightarrow0,
\]
so $\Ker(\eta_2)\simeq\Hom(\cU_X,\cQ_X^\vee)$. The spectral sequence degenerates at $E_2$ for degree reasons, and $\NHH^2(\cB_X,X)\simeq E_2^{-1,3}\simeq\Hom(\cU_X,\cQ_X^\vee)$.
\end{proof}

\begin{corollary}\label{cor:dimensions}
$\Ker(d\cP)$ has dimension $2$ for an ordinary and $3$ for a special Gushel--Mukai threefold, and dimension $5$ for a prime Fano threefold of genus eight.
\end{corollary}

\begin{proof}
By Proposition~\ref{prop:gammainj} and Theorem~\ref{thm:categorical}, $\Ker(d\cP)\simeq\Hom(\cU_X,\cQ_X^\vee)$, whose dimension is computed in \cite[Lemma~5.1]{JLLZ24} for Gushel--Mukai threefolds and in \cite[Lemma~A.8]{JLZ22a} in genus eight.
\end{proof}

For ordinary Gushel--Mukai threefolds, the two-dimensionality of $\Ker(d\cP)$ is \cite[Theorem~7.1]{DIM12}, which is used in the proof of \cite[Theorem~4.6]{JLLZ23}. For special Gushel--Mukai threefolds, Corollary~\ref{cor:dimensions} gives a second, categorical proof of Theorem~\ref{thm:kernelGM}. In both genera the computation works so smoothly because $\HH^1(\Ku(X))=0$: since $\HH^2(\Ku(X))$ is the tangent space to deformations of the category and $\HH^1(\Ku(X))$ its space of infinitesimal autoequivalences \cite[Sections~3.1--3.3]{Kuz15}, the vanishing of $\HH^1(\Ku(X))$ means that a normal Hochschild class, an infinitesimal change of the gluing between $\Ku(X)$ and the exceptional block, is already a genuine fixed-category deformation of $X$.

\begin{remark}\label{rem:othergenera}
For the remaining prime Fano threefolds of genus at least six, the categorical period map is also well understood. The map $\gamma_X$ is injective for $7\le g\le10$ by \cite[Theorem~4.10]{JLLZ23}, so again $\Ker(d\cP)=\Ker(\eta)$.
\begin{enumerate}[label=(\alph*)]
\item If $g=7,9,10$, then $\Ku(X)\simeq\Db(C)$ for a smooth curve $C$ of genus $7$, $3$, $2$ respectively \cite{Kuz06}, and the exact sequence of \cite[Theorem~8.8(2)]{Kuz09a}, together with the vanishing of $H^1(X,N^\vee_{X/\Gr}\otimes\omega_X^{-1})$, shows that $\eta\colon H^1(X,T_X)\to\HH^2(\Ku(X))\simeq H^1(C,T_C)$ is surjective. Since $h^1(X,T_X)=18$, $12$, $10$, the kernel has dimension $0$, $6$, $7$; in particular the infinitesimal Torelli theorem holds in genus seven.
\item If $g=12$, then $\Ku(X)$ is equivalent to the derived category of representations of the Kronecker quiver with three arrows \cite[Section~4]{Kuz09b}, which has $\HH^2=0$, the path algebra of an acyclic quiver being hereditary. Hence $\eta=0$ and $\Ker(d\cP)=H^1(X,T_X)\simeq\CC^6$; here $\gamma_X$ is trivially injective, since $H^1(X,\Omega^2_X)=H^2(X,\Omega^1_X)=0$.
\end{enumerate}
In these genera, however, $\HH^1(\Ku(X))\ne0$: it is $H^1(C,\cO_C)$, of dimension $7$, $3$, $2$ for $g=7,9,10$, and $\mathfrak{sl}_3$ for $g=12$. So normal Hochschild cohomology is strictly larger than $\Ker(\eta)$, and the argument of Theorem~\ref{thm:categorical} does not apply directly.
\end{remark}

\subsection{The kernel as a tangent space to a Bridgeland moduli space}
Let $X$ be a Gushel--Mukai threefold, $C(X)$ its Hilbert scheme of conics, and $\sigma$ a Serre-invariant stability condition on $\Ku(X)$. By \cite[Theorem~7.12]{JLLZ24}, the projection functor induces a birational morphism $C(X)\to M_\sigma(\Ku(X),[I_C])$, whose exceptional locus is a line or a plane according to whether $X$ is ordinary or special. We show that in both cases it is $\PP\Hom(\cU_X,\cQ_X^\vee)=\PP\Ker(d\cP)$, and that the same holds in genus eight.

\begin{theorem}\label{thm:Bridgeland}
Let $X$ be a prime Fano threefold of genus six or eight. Then $\PP\Ker(d\cP)\simeq\PP\Ker(\eta)\simeq\PP\Hom(\cU_X,\cQ_X^\vee)$ is identified with
\begin{enumerate}[label=\textup{(\arabic*)}]
\item the exceptional locus of the birational morphism $C(X)\to M_\sigma(\Ku(X),[I_C])$, if $g=6$;
\item the exceptional divisor of the birational morphism $M_X^{ss}(2,0,4)\to M_\sigma^{ss}(\Ku(X),2[I_C])$, where $M_X^{ss}(2,0,4)$ is the moduli space of semistable sheaves of rank two with $c_1=0$, $c_2=4$, $c_3=0$, if $g=8$.
\end{enumerate}
\end{theorem}

For genus eight we need the following results, obtained by Zhiyu Liu and the second author; for stability conditions on Kuznetsov components of Fano threefolds we refer to \cite[Section~3]{LZ22}.

\begin{lemma}\label{lem:Liu1}
Let $X$ be a prime Fano threefold of genus eight and $E\in\mathrm{Coh}^{-1/2}(X)$ a $\sigma_{\alpha,-1/2}$-semistable object for some $\alpha>0$ with $\ch(E)=(1,0,-2L,kP)$. Then $k\le0$.
\end{lemma}

\begin{proof}
The argument of \cite[Lemma~6.3]{LZ22} applies verbatim.
\end{proof}

\begin{lemma}\label{lem:Liu2}
Let $X$ be a prime Fano threefold of genus eight and $E\in M_X(2,0,4)$. Then $h^i(E)=0$ for all $i$, and $E\notin\Ku(X)$ if and only if there is a short exact sequence $0\to\cU_X\to\cQ_X^\vee\to E\to0$.
\end{lemma}

\begin{proof}
By \cite[Lemma~4.3]{BF11}, $h^2(E)=0$. Stability and Serre duality give $h^3(E)=\hom(E,\cO_X(-1))=0$ and $h^0(E)=0$, and $\chi(E)=0$ gives $h^1(E)=0$.

Suppose $E\notin\Ku(X)$. Stability and Serre duality give $\hom(\cU_X^\vee,E)=\mathrm{ext}^3(\cU_X^\vee,E)=0$, and $\chi(\cU_X^\vee,E)=0$, so $\mathrm{ext}^1(\cU_X^\vee,E)=\mathrm{ext}^2(\cU_X^\vee,E)\ne0$. Applying $\Hom(-,E)$ to $0\to\cQ_X^\vee\to\cO_X^{\oplus6}\to\cU_X^\vee\to0$ and using $h^i(E)=0$, we get $\Hom(\cQ_X^\vee,E)\ne0$; let $s\colon\cQ_X^\vee\to E$ be nonzero. Gieseker stability of $E$ gives $\ch_1(\im s)=0$, so $\ch_1(\ker s)=-H$, and $\ker s$ is $\mu$-stable by stability of $\cQ_X^\vee$.

If $\im s$ has rank one, then $\ch_2(\im s)=aL$ with $a\le-2$ and $\ch_3(\im s)=bP$ with $b\le-1$. By \cite[Proposition~4.8]{BBF+24}, $\im s$ and $E$ are $\sigma_{\alpha,-1/2}$-semistable of the same slope for $\alpha\gg0$, so $E/\im s$ is semistable; since $\ch_3(E/\im s)=-bP$ with $-b\ge1$, this contradicts Lemma~\ref{lem:Liu1} when $a=-2$. Hence $a\le-3$, so $\ch_{\le2}(\ker s)=(3,-H,(-a-2)L)$ with $-a-2\ge1$, contradicting \cite[Proposition~3.2]{Li19}. If $\im s$ has rank two, then $\ch_2(\im s)=aL$ with $a\le-4$, and $\ch_{\le2}(\ker s)=(2,-H,(-a-2)L)$; by \cite[Proposition~3.2]{Li19} we must have $a=-4$, so $\ch(\ker s)=(2,-H,2L,(c-\frac16)P)$ with $c\ge0$, and $\chi(\ker s)=0$ \cite[Proposition~3.5]{BF11} forces $c=0$. Thus $s$ is surjective and $\ker s\simeq\cU_X$.

Conversely, every nonzero $t\in\Hom(\cU_X,\cQ_X^\vee)\simeq\CC^5$ is injective by stability, and its cokernel $W$ is torsion-free and $\mu$-semistable. If $W\to W'$ were a minimal destabilizing quotient for Gieseker stability, then $\ch_{\le2}(W')=(1,0,xL)$ with $x\le-2$, and the kernel of $\cQ_X^\vee\to W'$ would have $\ch_{\le2}=(3,-H,(-x-2)L)$ with $-x-2\ge0$, again contradicting \cite[Proposition~3.2]{Li19}.
\end{proof}

\begin{proposition}\label{prop:Liu}
Let $X$ be a prime Fano threefold of genus eight and $\sigma$ a Serre-invariant stability condition on $\Ku(X)$. The projection functor induces a birational morphism $M_X^{ss}(2,0,4)\to M_\sigma^{ss}(\Ku(X),2[I_C])$ contracting a $\PP^4$ to a point.
\end{proposition}

\begin{proof}
The projection functor induces a morphism, which is projective since the source is projective and the target proper. It is \'etale and injective, hence an open immersion, on the open locus of sheaves lying in $\Ku(X)$, so it is birational. Since $\mathrm{ext}^1(\cQ_X^\vee,\cU_X[1])=1$, all sheaves $E\notin\Ku(X)$ are sent to a single point, and by Lemma~\ref{lem:Liu2} they are parametrized by $\PP\Hom(\cU_X,\cQ_X^\vee)\simeq\PP^4$.
\end{proof}

\begin{proof}[Proof of Theorem~\ref{thm:Bridgeland}]
(1) By \cite[Proposition~7.1]{JLLZ24}, $I_C\notin\Ku(X)$ exactly for the $\sigma$-conics of an ordinary $X$, parametrized by a line $L_\sigma$ \cite[p.~5]{DIM12}, and for the preimages of the lines of $Y$ when $X$ is special, parametrized by the plane $F(Y)\simeq\PP^2$, one of the two components of $C(X)$. For such conics there is an exact sequence $0\to\cU_X\to\cQ_X^\vee\to I_C\to0$. By \cite[Proposition~7.2 and Lemma~7.6]{JLLZ24}, the projection functor $L_{\cO_X}L_{\cU_X^\vee}$ contracts $C(X)$ along $\PP\Hom(\cU_X,\cQ_X^\vee)$, which is $L_\sigma\simeq\PP^1$ or $F(Y)\simeq\PP^2$. Theorem~\ref{thm:categorical} and Proposition~\ref{prop:gammainj} identify this with $\PP\Ker(d\cP)$.

(2) By Lemma~\ref{lem:Liu2}, the sheaves $E\notin\Ku(X)$ are the cokernels of $\cU_X\to\cQ_X^\vee$, parametrized by $\PP\Hom(\cU_X,\cQ_X^\vee)\simeq\PP^4$, and by Proposition~\ref{prop:Liu} this $\PP^4$ is the exceptional divisor of the contraction. Conclude as in (1).
\end{proof}

\begin{remark}\label{rem:notaccidental}
The identification is not accidental. For an ordinary Gushel--Mukai threefold, write $\Db(X)=\langle\Ku(X),\cQ_X^\vee,\cO_X\rangle$ with $\cQ_X^\vee=L_{\cO_X}\cU_X^\vee[-1]$, and consider the distinguished object $\pi(\cQ_X^\vee)=R_{\cU_X^\vee}\cQ_X^\vee\in\Ku(X)$. By \cite[Lemma~5.4]{JLLZ24} it is a two-term complex $\cU_X[1]\to\pi(\cQ_X^\vee)\to\cQ_X^\vee$, isomorphic to the projection of $I_C$ when $I_C\notin\Ku(X)$, and it is $\sigma$-stable \cite[Lemma~7.4]{JLLZ24}. Its tangent space is $\Ext^1(\pi(\cQ_X^\vee),\pi(\cQ_X^\vee))\simeq\Ext^1(\pi(\cQ_X^\vee),\cQ_X^\vee)$, and since $\cQ_X^\vee$ is exceptional, the long exact sequence obtained from the defining triangle gives $\Ext^1(\pi(\cQ_X^\vee),\pi(\cQ_X^\vee))\simeq\Hom(\cU_X,\cQ_X^\vee)$. Blowing up the moduli space at this point gives $C(X)$, and $\PP\Hom(\cU_X,\cQ_X^\vee)$ is the exceptional divisor. For special $X$ the point is singular, but the same identification of tangent spaces holds, and the genus-eight case is analogous.
\end{remark}

For a general ordinary Gushel--Mukai threefold, the connected component through $X$ of the fibre of the period map is $M_\sigma(\Ku(X),[I_C])$ modulo an involution \cite[Theorem~7.4]{DIM12}, \cite[Theorem~7.12]{JLLZ24}, and by Logachev's reconstruction theorem \cite[Theorem~7.7]{Log12} the point corresponding to $X$ is $[\pi(\cQ_X^\vee)]$. So $\Ker(d\cP)$ is the tangent space of this fibre at $X$, which recovers $\Ker(d\cP)\simeq\Hom(\cU_X,\cQ_X^\vee)$. For special Gushel--Mukai threefolds and in genus eight no Logachev-type theorem is available, but the distinguished object $E$ with $\cU_X[1]\to E\to\cQ_X^\vee$ still corresponds to the isomorphism class of $X$ by Brill--Noether reconstruction \cite[Theorems~1.1 and~1.3]{JLZ22a}.

\subsection{A direct answer to Question~\ref{q:direct}}
We now construct the isomorphism of Question~\ref{q:direct} without the period map, Hochschild cohomology or a description of the period fibre. Let $V$ be a five-dimensional vector space, $G=\Gr(2,V)$, with tautological sequence $0\to\cU_G\to V\otimes\cO_G\to\cQ_G\to0$, and let
\[
B\subset H^0(G,\cO_G(1))=\wedge^2V^\vee,\qquad\dim B=3,
\]
be the net of Pl\"ucker hyperplanes with $Y=G\cap\PP(B^\perp)$. Thus $N_{Y/G}^\vee\simeq B\otimes\cO_Y(-1)$. For the special Gushel--Mukai threefold $p\colon X\to Y$ we put $\cU_X=p^*\cU_Y$ and $\cQ_X=p^*\cQ_Y$.

\begin{theorem}\label{thm:polarity}
There is an $\Aut(Y)$-equivariant isomorphism, unique up to a nonzero scalar and independent of the branch surface,
\[
\Phi_X\colon H^1(Y,T_Y(-1))\xrightarrow{\ \sim\ }\Hom_X(\cU_X,\cQ_X^\vee).
\]
It is the composite of the natural isomorphisms $\delta\colon B^\vee\to H^1(Y,T_Y(-1))$ and $\beta_X\colon B\to\Hom_X(\cU_X,\cQ_X^\vee)$ defined below with the polarity $B\simeq B^\vee$ given by the conic of special lines on $Y$.
\end{theorem}

\begin{proof}
\emph{Step 1: $H^1(Y,T_Y(-1))\simeq B^\vee$.} The normal sequence twisted by $\cO_Y(-1)$ is
\begin{equation}\label{eq:normaltwist}
0\longrightarrow T_Y(-1)\longrightarrow T_G(-1)|_Y\longrightarrow B^\vee\otimes\cO_Y\longrightarrow0.
\end{equation}
We claim that $H^i(Y,T_G(-1)|_Y)=0$ for all $i$. The Koszul resolution of $\cO_Y$ expresses $T_G(-1)|_Y$ through $T_G(-1)$, $B\otimes T_G(-2)$, $\wedge^2B\otimes T_G(-3)$ and $\wedge^3B\otimes T_G(-4)$. In the standard convention the weight of $T_G(-j)=\cU_G^\vee\otimes\cQ_G\otimes\cO_G(-j)$ is $(1-j,-j\,|\,0,0,-1)$, and adding $\rho=(4,3\,|\,2,1,0)$ gives
\[
(4,2,2,1,-1),\quad(3,1,2,1,-1),\quad(2,0,2,1,-1),\quad(1,-1,2,1,-1)\qquad(j=1,2,3,4).
\]
Each has a repeated entry, so all four bundles are acyclic by the Borel--Weil--Bott theorem, which proves the claim. As $H^0(Y,\cO_Y)=\CC$ and $H^{>0}(Y,\cO_Y)=0$, the cohomology of \eqref{eq:normaltwist} gives the connecting isomorphism
\begin{equation}\label{eq:delta}
\delta\colon B^\vee\xrightarrow{\ \sim\ }H^1(Y,T_Y(-1)).
\end{equation}

\emph{Step 2: $B\simeq\Hom_Y(\cU_Y,\cQ_Y^\vee)$.} For $b\in B$ let $\sigma_b$ be the corresponding alternating form. At $[U]\in Y$ we have $b|_{\wedge^2U}=0$, so $u\mapsto(\bar v\mapsto\sigma_b(u,v))$ is a well-defined map $U\to(V/U)^\vee$; replacing $v$ by $v+u'$ changes the value by $\sigma_b(u,u')=0$. These maps define
\begin{equation}\label{eq:beta}
\beta_b\colon\cU_Y\longrightarrow\cQ_Y^\vee,\qquad\beta_Y\colon B\longrightarrow\Hom_Y(\cU_Y,\cQ_Y^\vee).
\end{equation}
To see that $\beta_Y$ is an isomorphism, put $F=\cU_G^\vee\otimes\cQ_G^\vee$. The shifted weights of $F(-j)$ for $j=0,1,2,3$ are
\[
(5,3,3,1,0),\quad(4,2,3,1,0),\quad(3,1,3,1,0),\quad(2,0,3,1,0).
\]
Those for $j=0,2,3$ are singular, and that for $j=1$ has one inversion and regularizes to $\rho$, so $H^1(G,F(-1))\simeq\CC$ and $H^i(G,F(-1))=0$ for $i\ne1$. In the Koszul spectral sequence for $F|_Y$ the only term of total degree zero is $B\otimes H^1(G,F(-1))$, hence
\begin{equation}\label{eq:HomB}
\Hom_Y(\cU_Y,\cQ_Y^\vee)=H^0(Y,F|_Y)\simeq B,
\end{equation}
and under this isomorphism $b$ corresponds to \eqref{eq:beta}; compare \cite[Lemma~3.3]{Fae05}.

\emph{Step 3: pulling back to $X$ gives no new morphisms.} Since the branch surface lies in $|\cO_Y(2)|$, $p_*\cO_X\simeq\cO_Y\oplus\cO_Y(-1)$, and the projection formula gives
\begin{equation}\label{eq:pushforward}
\Hom_X(\cU_X,\cQ_X^\vee)\simeq H^0(Y,F|_Y)\oplus H^0(Y,F(-1)|_Y).
\end{equation}
The Koszul resolution of $F(-1)|_Y$ involves $F(-1),\dots,F(-4)$; the first has cohomology only in degree one, and the others are acyclic, the shifted weight of $F(-4)$ being the singular weight $(1,-1,3,1,0)$. So the Koszul spectral sequence has no term of total degree zero, and $H^0(Y,F(-1)|_Y)=0$. Combining with \eqref{eq:HomB}, pullback of \eqref{eq:beta} gives an isomorphism
\begin{equation}\label{eq:betaX}
\beta_X\colon B\xrightarrow{\ \sim\ }\Hom_X(\cU_X,\cQ_X^\vee),
\end{equation}
which in particular does not depend on the branch surface.

\emph{Step 4: the conic of special lines identifies $B$ with $B^\vee$.} The automorphism group of $Y$ is $\PGL_2$, and $B$ is its irreducible three-dimensional representation; choosing a two-dimensional space $W$, we may write $B\simeq\Sym^2W$ \cite[Section~3]{Fae05}. The Hilbert scheme of lines on $Y$ is a projective plane, in which the special lines form a smooth conic \cite[Section~5.1]{KPS18}; under the equivariant identification with $\PP(B)$ it is the Veronese conic $\{[f^2]:[f]\in\PP(W)\}$. Its polarity is a nondegenerate $\PGL_2$-invariant symmetric form on $B$, unique up to scalar: for $q=ax^2+2bxy+cy^2$ and $q'=a'x^2+2b'xy+c'y^2$,
\begin{equation}\label{eq:kappa}
\kappa(q,q')=ac'+a'c-2bb',\qquad\kappa^\sharp\colon B\xrightarrow{\ \sim\ }B^\vee,\quad q\mapsto\kappa(q,-).
\end{equation}
We set $\Phi_X=\beta_X\circ(\kappa^\sharp)^{-1}\circ\delta^{-1}$. Every arrow is $\Aut(Y)$-equivariant, and since an invariant form on the irreducible representation $\Sym^2W$ is unique up to scalar, so is $\Phi_X$; its projectivization is canonical.
\end{proof}

\begin{remark}\label{rem:explicit}
Explicitly, given $\xi\in H^1(Y,T_Y(-1))$, let $\lambda_\xi=\delta^{-1}(\xi)\in B^\vee$ and let $b_\xi\in B$ be the unique element with $\lambda_\xi=\kappa(b_\xi,-)$, after normalizing $\kappa$. Then $\Phi_X(\xi)$ is the pullback to $X$ of $u\mapsto\sigma_{b_\xi}(u,-)|_{\cQ_Y}$. The construction converts a functional on the net of equations of $Y$ into the alternating form in that net polar to it.
\end{remark}

\begin{remark}\label{rem:family}
Let $\mathcal{Y}\to T$ be a family of embedded quintic del Pezzo threefolds, with rank-three bundle $\cB$ of Pl\"ucker equations. Relative versions of \eqref{eq:normaltwist} and \eqref{eq:beta}, with cohomology and base change, give $R^1\pi_*T_{\mathcal{Y}/T}(-1)\simeq\cB^\vee$ and $\pi_*\mathcal{H}om(\cU_{\mathcal{X}},\cQ_{\mathcal{X}}^\vee)\simeq\cB$. The relative conic of special lines determines a line subbundle $\mathcal{L}\subset\Sym^2\cB^\vee$ of nondegenerate forms, hence an isomorphism $\cB\to\cB^\vee\otimes\mathcal{L}^\vee$. So the construction is canonically functorial after projectivization, $\PP(R^1\pi_*T_{\mathcal{Y}/T}(-1))\simeq\PP(\pi_*\mathcal{H}om(\cU_{\mathcal{X}},\cQ_{\mathcal{X}}^\vee))$, while a linear isomorphism of vector bundles requires a trivialization of $\mathcal{L}$; with such a choice it commutes with arbitrary base change. This is the precise meaning of ``canonical up to scalar'' in Theorem~\ref{thm:polarity}.
\end{remark}

\begin{remark}\label{rem:twodescriptions}
For a special Gushel--Mukai threefold we now have
\[
H^1(Y,T_Y(-1))\simeq\Ker(d\cP)=\Ker(\eta)\simeq\Hom(\cU_X,\cQ_X^\vee).
\]
The left-hand space consists of infinitesimal deformations of $X$ in the direction normal to the special locus, towards possibly ordinary Gushel--Mukai threefolds with the same intermediate Jacobian. The right-hand space is the tangent space to the distinguished point $[\pi(\cQ_X^\vee)]$ of $M_\sigma(\Ku(X),[I_C])$, which is the connected component through $[X]$ of the fibre of the categorical period map \cite[Theorem~1.9]{BP23}, \cite[Theorem~1.1]{JLZ22a}; it parametrizes infinitesimal deformations of $X$ to period partners. Theorem~\ref{thm:polarity} relates the two descriptions directly: both come from the Pl\"ucker net $B$, and the conic of special lines exchanges the dual realizations $B^\vee$ and $B$.
\end{remark}

\section{Verra threefolds}\label{sec:Verra}

Verra threefolds are Fano threefolds of index one and degree twelve, related to Prym varieties of plane sextics \cite{Ver04}. Like Gushel--Mukai threefolds, they come in two kinds: \emph{ordinary} ones are smooth divisors of bidegree $(2,2)$ in $\PP^2\times\PP^2$, and \emph{special} ones are double covers of a divisor of bidegree $(1,1)$ branched along a K3 surface. In this section we study both kinds with the two methods developed above, and we interpret the kernel of the differential of the period map through the deformations of a spinor bundle.

The section is organized as follows. Section~\ref{ss:Verraprelim} introduces the tools common to both cases: a conic bundle structure, its Clifford component, the normal Hochschild cohomology of the exceptional collection, and spinor bundles. Sections~\ref{ss:Verraordinary} and~\ref{ss:Verraspecial} treat ordinary and special Verra threefolds separately, each time in the same order: first the infinitesimal period map, then the categorical period map, and finally the spinor bundle. Section~\ref{ss:Verracompare} compares the two cases.

\subsection{Preliminaries}\label{ss:Verraprelim}

\subsubsection{The two kinds of Verra threefolds}
Put $A=\PP^2\times\PP^2$.

\begin{definition}\label{def:Verra}
An \emph{ordinary Verra threefold} is a smooth divisor $X\in|\cO_A(2,2)|$. A \emph{special Verra threefold} is a double cover $\pi\colon X\to Y$ of a smooth divisor $Y=\{r=0\}\in|\cO_A(1,1)|$, branched along a smooth surface $S\in|\cO_Y(2,2)|=|-K_Y|$; the covering line bundle is $L=\cO_Y(1,1)$.
\end{definition}

In both cases we write $h$ and $k$ for the pullbacks of the two hyperplane classes of $A$, so that $H=h+k=-K_X$. In the special case $Y$ is the flag variety of $\PP^2$: it is a rigid Fano threefold with $\Pic(Y)=\mathbb{Z}^2$, $\omega_Y\simeq\cO_Y(-2,-2)$ and $H^1(Y,\Omega^2_Y)=0$, so the results of Section~\ref{sec:double} apply, and the branch surface $S$ is a K3 surface.

\subsubsection{The conic bundle structure and the Clifford component}
Let $f\colon X\to\PP^2$ be the projection to the first factor; the second projection can be treated in exactly the same way. In both cases $f$ is a flat conic bundle.

By Kuznetsov's description of the derived categories of quadric fibrations \cite[Theorem~4.2]{Kuz08}, there is a semiorthogonal decomposition
\begin{equation}\label{eq:VerraSOD}
\Db(X)=\langle\cC_f,\cO_X,\cO_X(h),\cO_X(2h)\rangle,\qquad\cC_f\simeq\Db(\PP^2,\cB_0),
\end{equation}
where $\cB_0$ is the sheaf of even Clifford algebras of the conic bundle, normalized so that the base block is $f^*\Db(\PP^2)$. We call $\cC_f$ the \emph{Clifford component}; all categories and functors carry their natural dg enhancements. We write $\cE=\langle\cO_X,\cO_X(h),\cO_X(2h)\rangle$ for the exceptional part. The Clifford component plays the role of the Kuznetsov component, and the corresponding infinitesimal categorical period map is
\[
\eta_f\colon H^1(X,T_X)\hookrightarrow\HH^2(X)\longrightarrow\HH^2(\cC_f),
\]
the composite of the inclusion of the geometric summand in the Hochschild--Kostant--Rosenberg decomposition $\HH^2(X)=H^2(X,\cO_X)\oplus H^1(X,T_X)\oplus H^0(X,\wedge^2T_X)$ with the restriction to $\cC_f$. It sends a first-order deformation of $X$ to the induced first-order deformation of $\cC_f$.

\subsubsection{Normal Hochschild cohomology of the exceptional collection}
The kernel of the restriction $\HH^2(X)\to\HH^2(\cC_f)$ is controlled by Kuznetsov's exact sequence \cite[Theorem~3.3]{Kuz15},
\begin{equation}\label{eq:VerraNHH}
\HH^1(\cC_f)\longrightarrow\NHH^2(\cE,X)\longrightarrow\HH^2(X)\longrightarrow\HH^2(\cC_f),
\end{equation}
so the first task is to compute the normal Hochschild cohomology $\NHH^2(\cE,X)$. Write
\[
V_h=H^0(X,\cO_X(h)),\qquad V_k=H^0(X,\cO_X(k)),\qquad V_H=H^0(X,\cO_X(H)),
\]
and let $\mu\colon V_h\otimes V_k\to V_H$ be the multiplication map. We need some vanishing results.

\begin{lemma}\label{lem:Verravanish}
For ordinary and for special Verra threefolds, the exceptional collection $\cE$ is strong, $H^{>0}(X,\cO_X(H))=H^{>0}(X,\cO_X(k))=0$, and $H^\bullet(X,\cO_X(-h+k))=0$.
\end{lemma}

\begin{proof}
For an ordinary Verra threefold, these statements follow from the restriction sequences $0\to\cO_A(D-(2,2))\to\cO_A(D)\to\cO_X(D)\to0$ and the cohomology of line bundles on $A$. For example, the sequence $0\to\cO_A(-3,-1)\to\cO_A(-1,1)\to\cO_X(-h+k)\to0$ has acyclic outer terms. For a special Verra threefold, the formula $\pi_*\cO_X=\cO_Y\oplus\cO_Y(-1,-1)$ reduces the claims to line bundles on $Y$, and the restriction sequences of $Y\subset A$ reduce them further to line bundles on $A$. For example, $\cO_X(-h+k)$ gives the line bundles $\cO_Y(-1,1)$ and $\cO_Y(-2,0)$ on $Y$; their restriction sequences have acyclic left terms $\cO_A(-2,0)$ and $\cO_A(-3,-1)$ and acyclic middle terms.
\end{proof}

\begin{proposition}\label{prop:VerraNHH}
For ordinary and for special Verra threefolds, $\NHH^2(\cE,X)\simeq(\Ker\mu)^{\oplus2}$.
\end{proposition}

\begin{proof}
Put $E_i=\cO_X(ih)$ for $i=0,1,2$. The normal Hochschild spectral sequence \cite[Proposition~3.7]{Kuz15} has first page
\begin{multline}\label{eq:VerraE1}
E_1^{-p,q}=\bigoplus_{0\le i_0<\dots<i_p\le2}\ \bigoplus_{k_0+\dots+k_p=q}\Ext^{k_0}(E_{i_0},E_{i_1})\otimes\cdots\\
\cdots\otimes\Ext^{k_{p-1}}(E_{i_{p-1}},E_{i_p})\otimes H^{k_p-3}\bigl(X,\cO_X(H+(i_0-i_p)h)\bigr),
\end{multline}
and converges to normal Hochschild cohomology in degree $q-p$; the last factor comes from the inverse Serre functor $S_X^{-1}=(-)\otimes\omega_X^{-1}[-3]$. By Lemma~\ref{lem:Verravanish}, all terms vanish except in the row $q=3$, and there only the columns $p=0,1$ survive:
\[
E_1^{-1,3}=(V_h\otimes V_k)^{\oplus2},\qquad E_1^{0,3}=V_H^{\oplus3}.
\]
Here the two summands of $E_1^{-1,3}$ come from the pairs $(0,1)$ and $(1,2)$, while the pair $(0,2)$ and the triple $(0,1,2)$ contribute nothing, since their last factor is the cohomology of $\cO_X(-h+k)$. Up to the signs of the bar complex, the only differential is
\[
d_1(a,b)=\bigl(-\mu(a),\ \mu(a)-\mu(b),\ \mu(b)\bigr),
\]
whose kernel is $(\Ker\mu)^{\oplus2}$, and there are no other differentials.
\end{proof}

\subsubsection{Spinor bundles}
Following Kuznetsov \cite[Definition~1.1]{Kuz25}, an \emph{abstract spinor bundle} for the conic bundle $f$ is a vector bundle $\cF$ of rank two on $X$ such that
\[
Rf_*\cF=0\qquad\text{and}\qquad c_1(\cF)=K_{X/\PP^2}\ \text{ in }\Pic(X/\PP^2).
\]
We use the \emph{canonical spinor bundle}, the middle term of the unique nonsplit extension
\begin{equation}\label{eq:spinor}
0\longrightarrow\omega_{X/\PP^2}\longrightarrow\cF\longrightarrow\cO_X\longrightarrow0
\end{equation}
with vanishing derived pushforward \cite[Lemma~2.16]{Kuz25}. It differs from the spinor bundle attached to the regular Clifford module of an unnormalized quadratic presentation by a line bundle pulled back from the base, which changes neither its endomorphism algebra nor its deformations. The notion of a spinor bundle depends on the chosen projection $f$.

Let $\mathrm{Spin}_f(X)$ be the local moduli stack of abstract spinor bundles on $f$ with the Chern classes of $\cF$, that is, the substack of the stack of rank-two vector bundles on $X$, near $[\cF]$, cut out by the spinor conditions. The deformation theory of a vector bundle is governed by its self-extensions: first-order deformations of $\cF$ are classified by $\Ext^1_X(\cF,\cF)$, and obstructions lie in $\Ext^2_X(\cF,\cF)$ \cite[Section~7]{Har10}, \cite[Appendix~2.A]{HL10}. After rigidifying the scalar automorphisms, the space of first-order deformations of $\mathrm{Spin}_f(X)$ at $[\cF]$ is therefore
\begin{equation}\label{eq:Spintangent}
T_{[\cF]}\mathrm{Spin}_f(X)=\Ext^1_X(\cF,\cF).
\end{equation}
Here the spinor conditions impose no further condition to first order: a deformation of a vector bundle remains locally free, its derived pushforward remains zero by the exact sequence over the dual numbers, and its determinant cannot move because $H^1(X,\cO_X)=0$. We make no claim that $\mathrm{Spin}_f(X)$ is smooth; in the special case it is not.

The self-extensions of $\cF$ can be computed on the base. For a quadratic presentation $q\colon\mathcal{L}\to\Sym^2\cE_0^\vee$ of the conic bundle, the even Clifford algebra is $\cB_0=\cO_{\PP^2}\oplus(\wedge^2\cE_0\otimes\mathcal{L})$, and the Clifford equivalence identifies the canonical spinor bundle with the regular Clifford module, up to a twist from the base. Hence
\begin{equation}\label{eq:Cliffordext}
Rf_*\mathcal{E}nd(\cF)\simeq\cB_0,\qquad\Ext^i_X(\cF,\cF)\simeq H^i(\PP^2,\cB_0),
\end{equation}
compatibly with multiplication, since the first isomorphism is induced by the Clifford equivalence and by composition of endomorphisms; compare \cite[Corollary~2.11 and equation~(16)]{Kuz25}.

\subsection{Ordinary Verra threefolds}\label{ss:Verraordinary}
Let $X\subset A$ be an ordinary Verra threefold. Iliev proved the global Torelli theorem for a general Verra threefold \cite[Theorem~5.6]{Ili97}. We show that the infinitesimal Torelli theorem holds for every ordinary Verra threefold, that the categorical period map is injective, and that the canonical spinor bundle is rigid.

\subsubsection{The infinitesimal Torelli theorem}
We follow Flenner's method \cite[Theorem~1.1]{Fle86}, as in \cite[Theorem~5.1]{DIM12}. The idea is to use connecting maps to move the cup product $H^1(X,T_X)\otimes H^1(X,\Omega^2_X)\to H^2(X,\Omega^1_X)$ into the top cohomology of line bundles, where it becomes Serre dual to a multiplication map of sections. The conormal sequence of $X\subset A$ and \cite[Theorem~2.1]{Fle86} give the three exact sequences
\begin{align}
&0\longrightarrow\cO_X(-2,-2)\longrightarrow\Omega^1_A|_X\longrightarrow\Omega^1_X\longrightarrow0,\label{eq:Vconormal}\\
&0\longrightarrow\Omega^1_X(-1,-1)\longrightarrow\Omega^2_A(1,1)|_X\longrightarrow T_X\longrightarrow0,\label{eq:VFlenner}\\
&0\longrightarrow\Omega^2_X\longrightarrow(\Omega^1_A|_X)^\vee\otimes\omega_X\longrightarrow\cO_X(1,1)\longrightarrow0.\label{eq:Vdual}
\end{align}
The second is obtained by twisting the sequence $0\to\Omega^1_X(-2,-2)\to\Omega^2_A|_X\to\Omega^2_X\to0$ by $\omega_X^{-1}\simeq\cO_X(1,1)$ and using $T_X\otimes\omega_X\simeq\Omega^2_X$; the third is the dual of \eqref{eq:Vconormal} twisted by $\omega_X$. Denote by
\begin{align*}
\partial_2&\colon H^1(X,T_X)\longrightarrow H^2\bigl(X,\Omega^1_X(-1,-1)\bigr),\\
\partial_1&\colon H^2\bigl(X,\Omega^1_X(-1,-1)\bigr)\longrightarrow H^3\bigl(X,\cO_X(-3,-3)\bigr),\\
\delta&\colon H^0\bigl(X,\cO_X(1,1)\bigr)\longrightarrow H^1(X,\Omega^2_X)
\end{align*}
the connecting maps of \eqref{eq:VFlenner}, of \eqref{eq:Vconormal} twisted by $\cO_X(-1,-1)$, and of \eqref{eq:Vdual}. By \cite[Lemma~2.10]{Fle86}, these maps are compatible with cup products up to sign; in particular, for $\tau\in H^1(X,T_X)$ and $\gamma\in H^0(X,\cO_X(1,1))$, the image of $\tau\cdot\delta(\gamma)\in H^2(X,\Omega^1_X)$ under the connecting map of \eqref{eq:Vconormal} is, up to sign,
\begin{equation}\label{eq:Vproduct}
\partial_1\partial_2(\tau)\cdot\gamma\ \in H^3\bigl(X,\cO_X(-2,-2)\bigr).
\end{equation}
Hence $d\cP$ is injective as soon as the following two conditions hold:
\begin{enumerate}[label=(\alph*)]
\item the maps $\partial_2$ and $\partial_1$ are injective;
\item the pairing $H^3(X,\cO_X(-3,-3))\otimes H^0(X,\cO_X(1,1))\to H^3(X,\cO_X(-2,-2))$ is nondegenerate in the first factor.
\end{enumerate}
Indeed, if $\tau\ne0$, then (a) gives $\partial_1\partial_2(\tau)\ne0$, (b) gives a $\gamma$ for which \eqref{eq:Vproduct} is nonzero, and therefore $\tau\cdot\delta(\gamma)\ne0$.

Condition (b) is Serre dual to the surjectivity of the multiplication map $H^0(X,\cO_X(1,1))\otimes H^0(X,\cO_X(1,1))\to H^0(X,\cO_X(2,2))$. This holds because the restriction maps from $A$ are surjective and multiplication is surjective on $A$ by the K\"unneth formula \cite[p.~317]{Fle86}. Condition (a) follows from the long exact sequences of \eqref{eq:VFlenner} and of \eqref{eq:Vconormal} twisted by $\cO_X(-1,-1)$, together with the following vanishings.

\begin{lemma}\label{lem:Vvanish}
For an ordinary Verra threefold,
\[
H^1\bigl(X,\Omega^2_A(1,1)|_X\bigr)=0\qquad\text{and}\qquad H^2\bigl(X,\Omega^1_A(-1,-1)|_X\bigr)=0.
\]
\end{lemma}

\begin{proof}
By the ideal sequence of $X\subset A$, the first vanishing follows from $H^1(A,\Omega^2_A(1,1))=0$ and $H^2(A,\Omega^2_A(-1,-1))=0$, and the second from $H^2(A,\Omega^1_A(-1,-1))=0$ and $H^3(A,\Omega^1_A(-3,-3))=0$. All four follow from the K\"unneth formula, using for instance
\[
\Omega^2_A(1,1)\simeq\bigl(\Omega^2_{\PP^2}(1)\boxtimes\cO(1)\bigr)\oplus\bigl(\Omega^1_{\PP^2}(1)\boxtimes\Omega^1_{\PP^2}(1)\bigr)\oplus\bigl(\cO(1)\boxtimes\Omega^2_{\PP^2}(1)\bigr),
\]
and from Kodaira--Akizuki--Nakano vanishing.
\end{proof}

We now come to the first main result of this section.

\begin{proposition}\label{prop:VerraTorelli}
For every ordinary Verra threefold $X$, the infinitesimal period map
\[
d\cP\colon H^1(X,T_X)\longrightarrow\Hom\bigl(H^1(X,\Omega^2_X),H^2(X,\Omega^1_X)\bigr)
\]
is injective.
\end{proposition}

\begin{proof}
Conditions (a) and (b) hold by Lemma~\ref{lem:Vvanish} and the surjectivity of multiplication, so $d\cP$ is injective by the argument following \eqref{eq:Vproduct}.
\end{proof}

\subsubsection{The categorical period map}
{\sloppy For an ordinary Verra threefold the restriction map from $A$ identifies $V_H$ with $H^0(A,\cO_A(1,1))=V_h\otimes V_k$, so $\mu$ is an isomorphism, and Proposition~\ref{prop:VerraNHH} gives $\NHH^2(\cE,X)=0$. The exact sequence \eqref{eq:VerraNHH} then shows that the restriction $\HH^2(X)\to\HH^2(\cC_f)$ is injective on all of $\HH^2(X)$, and in particular on its geometric summand.\par}

\begin{proposition}\label{prop:Verraeta}
For an ordinary Verra threefold, $\eta_f$ is injective, and $X$ has no nonzero global bivector fields: $H^0(X,\wedge^2T_X)=0$.
\end{proposition}

\begin{proof}
Injectivity of $\eta_f$ was just shown. For the second assertion, note that $\wedge^2T_X\simeq\Omega^1_X(H)$, since $-K_X=H$. In the twisted conormal sequence
\[
0\longrightarrow\cO_X(-H)\longrightarrow\Omega^1_A(H)|_X\longrightarrow\Omega^1_X(H)\longrightarrow0,
\]
the middle term is acyclic, because the two terms $\Omega^1_A(-1,-1)$ and $\Omega^1_A(1,1)$ of its restriction sequence are acyclic by the Euler sequences; and $H^1(X,\cO_X(-H))=0$ by the restriction sequence with terms $\cO_A(-3,-3)$ and $\cO_A(-1,-1)$. Hence $H^0(X,\Omega^1_X(H))=0$.
\end{proof}

\subsubsection{Rigidity of the canonical spinor bundle}
For an ordinary Verra threefold, the projection $f$ presents $X$ as a relative quadric in $\PP^2\times\PP^2$ over the first factor, with
\[
\cE_0=\cO_{\PP^2}^{\oplus3},\qquad\mathcal{L}=\cO_{\PP^2}(-2),\qquad\cB_0=\cO_{\PP^2}\oplus\cO_{\PP^2}(-2)^{\oplus3}.
\]
Since $\cO_{\PP^2}(-2)$ has no cohomology, \eqref{eq:Cliffordext} gives
\begin{equation}\label{eq:ordinaryExt}
\Ext^0_X(\cF,\cF)=\CC,\qquad\Ext^i_X(\cF,\cF)=0\quad(i>0).
\end{equation}
So the canonical spinor bundle is simple and rigid, and it is an isolated reduced point of its local moduli space. Combining this with Propositions~\ref{prop:VerraTorelli} and~\ref{prop:Verraeta}, all three points of view give the same, trivial, answer:
\begin{equation}\label{eq:ordinarysummary}
T_{[\cF]}\mathrm{Spin}_f(X)=0=\Ker(\eta_f)=\Ker(d\cP).
\end{equation}

\begin{remark}\label{rem:geometricordinary}
The agreement in \eqref{eq:ordinarysummary} has a geometric explanation. The kernel of $\eta_f$ is the tangent space to the fibre of the categorical period map through $[X]$, that is, the space of first-order deformations of $X$ which leave the Clifford component $\cC_f$ unchanged. By Kuznetsov's modification theorem \cite[Theorem~1.2]{Kuz25}, conic bundles over $\PP^2$ with the same Clifford category correspond to abstract spinor bundles, so such deformations of $X$ correspond to deformations of $\cF$, whose tangent space is $\Ext^1_X(\cF,\cF)$; Lemma~\ref{lem:spinormod} below makes this correspondence explicit at first order. Since $\cF$ is rigid, there are no such deformations, which is the vanishing of $\Ker(\eta_f)$.

Two points deserve emphasis. First, $\Ker(d\cP)$ is the tangent space to the \emph{scheme-theoretic} fibre of the period map, and a fibre consisting of a single point may still have a nonzero tangent space; this is what happens for hyperelliptic curves, and for special Verra threefolds, where the local moduli of $\cF$ is $\Spec\CC[t]/(t^2)$ (Proposition~\ref{prop:obstruction}). So the vanishing of $\Ker(d\cP)$ cannot be deduced from the global Torelli theorem of Iliev \cite[Theorem~5.6]{Ili97}, which moreover concerns general Verra threefolds only; Proposition~\ref{prop:VerraTorelli} proves it directly for every ordinary Verra threefold. Second, in general only the inclusion $\Ker(\eta_f)\subseteq\Ker(d\cP)$ follows from the factorization of $d\cP$ through Hochschild cohomology; the equality in \eqref{eq:ordinarysummary} uses the two independent computations.
\end{remark}

\subsection{Special Verra threefolds}\label{ss:Verraspecial}
Let $\pi\colon X\to Y$ be a special Verra threefold with covering involution $\iota$, and let $p\colon Y\to\PP^2$ be the first projection, so that $f=p\circ\pi$. In contrast with the ordinary case, the involution now forces a one-dimensional kernel. We determine this kernel by the Hodge-theoretic method of Section~\ref{sec:GM}, compute it again categorically, and identify it with the deformations of the canonical spinor bundle, which turn out to be obstructed.

\subsubsection{The decomposition of the deformation space}

\begin{lemma}\label{lem:Vdecomp}
The invariant and anti-invariant parts of $H^1(X,T_X)$ have dimensions
\[
h^1\bigl(Y,T_Y(-\log S)\bigr)=18,\qquad h^1\bigl(Y,T_Y(-1,-1)\bigr)=1,
\]
so $h^1(X,T_X)=19$. The anti-invariant part $H^1(X,T_X)^-=H^1(Y,T_Y(-1,-1))$ lies in $\Ker(d\cP)$.
\end{lemma}

\begin{proof}
As in Lemma~\ref{lem:antiinv}, $H^2(Y,T_Y(-\log S))=0$, since $H^1(S,\cO_S(2,2))=0$ and $H^2(Y,T_Y)=0$. Since $T_Y(-S)\simeq\Omega^2_Y$, the sequence \eqref{eq:logTS} gives
\[
0\longrightarrow H^1\bigl(Y,T_Y(-\log S)\bigr)\longrightarrow H^1(S,T_S)\longrightarrow H^2(Y,\Omega^2_Y)\longrightarrow0,
\]
and $h^{2,2}(Y)=2$, so $h^1(Y,T_Y(-\log S))=20-2=18$. For the anti-invariant part, consider the normal sequence of $Y\subset A$ twisted by $\cO_Y(-1,-1)$,
\[
0\longrightarrow T_Y(-1,-1)\longrightarrow T_A(-1,-1)|_Y\longrightarrow\cO_Y\longrightarrow0.
\]
Its middle term is acyclic, since both terms $T_A(-1,-1)$ and $T_A(-2,-2)$ of its restriction sequence are acyclic by the Euler sequences and the K\"unneth formula; so the connecting map gives $H^0(Y,\cO_Y)\simeq H^1(Y,T_Y(-1,-1))$. Finally, $H^1(Y,\Omega^2_Y)=0$, so the anti-invariant part lies in the kernel by Remark~\ref{rem:parity}.
\end{proof}

\subsubsection{The kernel of the infinitesimal period map}
It remains to show that the invariant part of $d\cP$ is injective. We follow the three steps of Section~\ref{sec:GM}: we show that the relevant residue map is surjective, reduce the problem to a twisted infinitesimal Torelli statement for the branch K3 surface $S$, and prove that statement by Flenner's method. The only new feature is that the base now has Picard rank two.

\emph{Step 1: surjectivity of the residue.} By Lemma~\ref{lem:sign} twisted by $\cO_Y(-1,-1)$, the invariant part of $d\cP$ is compatible, up to sign, with the cup product on $S$, through the injection $H^1(Y,T_Y(-\log S))\hookrightarrow H^1(S,T_S)$ with image $H^1(S,T_S)_0$ and the residue
\[
\res\colon H^1\bigl(Y,\Omega^2_Y(\log S)(-1,-1)\bigr)\longrightarrow H^1\bigl(S,\Omega^1_S(-1,-1)\bigr).
\]
This residue is surjective because $H^2(Y,\Omega^2_Y(-1,-1))=0$. To see this, note that $\Omega^2_Y(-1,-1)\simeq T_Y(-3,-3)$ and use the normal sequence $0\to T_Y(-3,-3)\to T_A(-3,-3)|_Y\to\cO_Y(-2,-2)\to0$. Since $H^1(Y,\cO_Y(-2,-2))=0$, it suffices to show $H^2(Y,T_A(-3,-3)|_Y)=0$, which is Serre dual to $H^1(Y,\Omega^1_A(1,1)|_Y)=0$. The latter follows from the sequence $0\to\Omega^1_A\to\Omega^1_A(1,1)\to\Omega^1_A(1,1)|_Y\to0$, since $H^2(A,\Omega^1_A)=0$ and $H^1(A,\Omega^1_A(1,1))=0$ by the K\"unneth formula applied to $\Omega^1_A(1,1)\simeq(\Omega^1_{\PP^2}(1)\boxtimes\cO(1))\oplus(\cO(1)\boxtimes\Omega^1_{\PP^2}(1))$. As in Lemma~\ref{lem:reduction}, the injectivity of the invariant part of $d\cP$ is thus reduced to the nondegeneracy, in the first factor, of the pairing
\begin{equation}\label{eq:Vpairing}
H^1(S,T_S)_0\otimes H^1\bigl(S,\Omega^1_S(-1,-1)\bigr)\longrightarrow H^2\bigl(S,\cO_S(-1,-1)\bigr).
\end{equation}

\emph{Step 2: Flenner's method on $S$.} Since $N_{S/Y}\simeq\cO_S(2,2)$, the conormal sequence of $S\subset Y$ and \cite[Theorem~2.1]{Fle86} give
\[
0\to\cO_S(-2,-2)\to\Omega^1_Y|_S\to\Omega^1_S\to0,\qquad0\to\Omega^1_S(-1,-1)\to\Omega^2_Y(1,1)|_S\to\Omega^2_S(1,1)\to0.
\]
Fix a nowhere vanishing $\omega\in H^0(S,\omega_S)$ and let
\[
I\colon H^1(S,T_S)_0\hookrightarrow H^1(S,T_S)\xrightarrow{\ \lrcorner\omega\ }H^1(S,\Omega^1_S)\longrightarrow H^2\bigl(S,\cO_S(-2,-2)\bigr),
\]
the last map being the connecting map of the first sequence. As in Lemma~\ref{lem:Flennerdiagram}, \cite[Lemma~2.10]{Fle86} shows that \eqref{eq:Vpairing} is nondegenerate in the first factor as soon as $I$ is injective and the multiplication pairing
\[
H^2\bigl(S,\cO_S(-2,-2)\bigr)\otimes H^0\bigl(S,\cO_S(1,1)\bigr)\longrightarrow H^2\bigl(S,\cO_S(-1,-1)\bigr)
\]
is nondegenerate in the first factor. The latter is Serre dual to the surjectivity of $H^0(S,\cO_S(1,1))^{\otimes2}\to H^0(S,\cO_S(2,2))$, which, as in Lemma~\ref{lem:Noether}, reduces to Noether's theorem on a general member of $|\cO_S(1,1)|$, a smooth canonical curve.

\emph{Step 3: injectivity of $I$.} The kernel of $H^1(S,\Omega^1_S)\to H^2(S,\cO_S(-2,-2))$ is the image of $J_2\colon H^1(S,\Omega^1_Y|_S)\to H^1(S,\Omega^1_S)$, which is injective because $H^1(S,\cO_S(-2,-2))=0$. As in Lemma~\ref{lem:H1restr}, $H^1(S,\Omega^1_Y|_S)\simeq H^1(Y,\Omega^1_Y)\simeq\CC^2$, so this image is spanned by the restrictions $h_S$ and $k_S$ of the two hyperplane classes. Suppose that $0\ne\tau\in\im\mathrm{KS}$ satisfies $\tau\lrcorner\omega=J_2(a)$. Since Kodaira--Spencer classes map to zero in $H^1(S,T_Y|_S)$, the argument of Lemma~\ref{lem:Iinj} gives $\int_SJ_2(a)\wedge J_2(b)=0$ for every $b$. But the intersection form on the span of $h_S$ and $k_S$ has Gram matrix
\[
\begin{pmatrix}h_S^2&h_Sk_S\\h_Sk_S&k_S^2\end{pmatrix}=\begin{pmatrix}2&4\\4&2\end{pmatrix},
\]
which is nondegenerate; hence $J_2(a)=0$ and $\tau=0$. On all of $H^1(S,T_S)$ the corresponding map has a two-dimensional kernel, so, as for Gushel--Mukai threefolds, the restriction to $H^1(S,T_S)_0$ is essential.

We can now state the main Hodge-theoretic result for special Verra threefolds.

\begin{proposition}\label{prop:specialVerra}
For a special Verra threefold, the invariant infinitesimal period map
\[
H^1\bigl(Y,T_Y(-\log S)\bigr)\longrightarrow\Hom\Bigl(H^1\bigl(Y,\Omega^2_Y(\log S)(-1,-1)\bigr),H^2\bigl(Y,\Omega^1_Y(\log S)(-1,-1)\bigr)\Bigr)
\]
is injective. Consequently $\Ker(d\cP)=H^1(Y,T_Y(-1,-1))$, a one-dimensional space.
\end{proposition}

\begin{proof}
Steps 1--3 prove the injectivity of the invariant part, and Lemma~\ref{lem:Vdecomp} identifies the kernel with the anti-invariant part.
\end{proof}

\subsubsection{The categorical period map}
We now compute the same kernel categorically. Two features distinguish the special case from the ordinary one. First, the exceptional collection now has nonzero normal Hochschild cohomology. By the projection formula, $V_H=H^0(Y,\cO_Y(1,1))\oplus H^0(Y,\cO_Y)$, of dimension $8+1=9$; the image of $\mu$ is the first summand, the second summand is spanned by the coordinate of the double cover, and $\Ker\mu=\CC r$ is spanned by the equation of $Y$. Hence $\NHH^2(\cE,X)\simeq\CC^2$ by Proposition~\ref{prop:VerraNHH}. Second, $X$ carries a global bivector field, so the full group $\HH^2(X)$ is larger than its geometric summand.

\begin{lemma}\label{lem:bivector}
For a special Verra threefold, $H^0(X,\wedge^2T_X)$ is one-dimensional and anti-invariant under $\iota$.
\end{lemma}

\begin{proof}
As $-K_X=H$, we have $\wedge^2T_X\simeq\Omega^1_X(H)$, and the double cover formula gives
\[
H^0\bigl(X,\Omega^1_X\otimes\pi^*\cO_Y(1,1)\bigr)=H^0\bigl(Y,\Omega^1_Y(1,1)\bigr)\oplus H^0\bigl(Y,\Omega^1_Y(\log S)\bigr).
\]
The second term vanishes, since the boundary map of the residue sequence sends $1\in H^0(S,\cO_S)$ to $c_1(S)\ne0$. For the first term, the sequences $0\to\Omega^1_A\to\Omega^1_A(1,1)\to\Omega^1_A(1,1)|_Y\to0$ and $0\to\cO_Y\to\Omega^1_A(1,1)|_Y\to\Omega^1_Y(1,1)\to0$ give $H^0(Y,\Omega^1_A(1,1)|_Y)\simeq H^1(A,\Omega^1_A)\simeq\CC^2$ and $h^0(Y,\Omega^1_Y(1,1))=1$. To determine the action of $\iota$, note that the natural linearization of $\omega_X^{-1}$ differs from the pullback linearization of $\pi^*\cO_Y(1,1)$ by the sign character: at a ramification point, $\iota$ acts on $T_X$ with eigenvalues $(-1,1,1)$, but trivially on the fibre of a pullback line bundle. The nonzero section is invariant for the pullback linearization, so the bivector field is anti-invariant.
\end{proof}

The key input is the vanishing of $\HH^1(\cC_f)$, which follows from a symmetry of the Clifford component: $\cC_f$ is an \emph{Enriques category}, that is, its Serre functor is $\iota[2]$, and its category of $\iota$-invariant objects is $\Db(S)$ \cite[Example~3.11]{PPZ26}, \cite[Theorem~1.1 and Proposition~7.10]{KP17}. We now come to the categorical counterpart of Proposition~\ref{prop:specialVerra}.

\begin{theorem}\label{thm:Verracategorical}
For a special Verra threefold,
\begin{equation}\label{eq:Verrakernel}
\Ker(\eta_f)=H^1(X,T_X)^-=H^1(Y,T_Y(-1,-1))\simeq\CC.
\end{equation}
Moreover, the kernel of the full restriction $\HH^2(X)\to\HH^2(\cC_f)$ is the two-dimensional space
\begin{equation}\label{eq:fullkernel}
\HH^2(X)^-=H^1(X,T_X)^-\oplus H^0(X,\wedge^2T_X).
\end{equation}
\end{theorem}

\begin{proof}
We first show that $\HH^1(\cC_f)=0$. By \cite[Theorem~1.1]{Per21}, $\HH^\bullet(\cC_f)^{\langle\iota\rangle}$ is a direct summand of $\HH^\bullet(\Db(S))$. The action of $\iota$ on $\HH^\bullet(\cC_f)$ is trivial, since $\iota$ is the Serre functor up to an even shift and Serre functors act trivially on Hochschild cohomology \cite[Proposition~1.2]{Per21}. As $\HH^1(S)=0$, we get $\HH^1(\cC_f)=0$, so the map $\NHH^2(\cE,X)\to\HH^2(X)$ in \eqref{eq:VerraNHH} is injective.

Next, the two normal Hochschild classes are anti-invariant. All forward arrows in the bar complex carry pullback linearizations, each return arrow carries one factor of the inverse Serre functor, and the return sections are sections of $\cO(k)$ pulled back from $Y$; so the sign change of Lemma~\ref{lem:bivector} makes these classes anti-invariant. Since \eqref{eq:VerraNHH} is equivariant, the image of $\NHH^2(\cE,X)$ lies in $\HH^2(X)^-$.

Finally, $H^2(X,\cO_X)=0$, so Lemmas~\ref{lem:Vdecomp} and~\ref{lem:bivector} show that $\HH^2(X)^-$ is the two-dimensional space \eqref{eq:fullkernel}. The injective map from the two-dimensional space $\NHH^2(\cE,X)$ therefore has image exactly $\HH^2(X)^-$, which is the kernel of the restriction to $\cC_f$. Intersecting with the geometric summand $H^1(X,T_X)$ gives \eqref{eq:Verrakernel}.
\end{proof}

Comparing with Proposition~\ref{prop:specialVerra}, we obtain $\Ker(\eta_f)=\Ker(d\cP)$. The two computations are independent, but the equality uses both of them: the factorization of $d\cP$ through the Hochschild action only gives the inclusion $\Ker(\eta_f)\subseteq\Ker(d\cP)$. Note also that normal Hochschild cohomology alone does not compute the kernel of the differential of the period map here: it contains the class of the bivector field as well, and one has to intersect with the geometric summand.

\subsubsection{Deformations of the canonical spinor bundle}
We now interpret the one-dimensional kernel of the differential of the period map in terms of the canonical spinor bundle. Identify $Y$ with $\PP_{\PP^2}(\Omega^1_{\PP^2})$, whose points are lines in the fibres, so that $\cO_{\PP(\Omega^1)}(1)\simeq\cO_Y(1,1)$. The double cover branched in $|\cO_Y(2,2)|$ is then the relative quadric in $\PP_{\PP^2}(\Omega^1_{\PP^2}\oplus\cO_{\PP^2})$ given by
\begin{equation}\label{eq:relquadric}
q(u,z)=b(u,u)-z^2=0,
\end{equation}
where $b\colon\Sym^2\Omega^1_{\PP^2}\to\cO_{\PP^2}$ is the symmetric form defining the branch divisor. Thus $\cE_0=\Omega^1_{\PP^2}\oplus\cO_{\PP^2}$, $\mathcal{L}=\cO_{\PP^2}$, and
\begin{equation}\label{eq:specialB0}
\cB_0=\cO_{\PP^2}\oplus\Omega^1_{\PP^2}\oplus\Omega^2_{\PP^2}.
\end{equation}
The covering involution is $z\mapsto-z$; it acts by $-1$ on the summand $\Omega^1_{\PP^2}$ and trivially on the other two. By \eqref{eq:Cliffordext} and the cohomology of differential forms on $\PP^2$, the canonical spinor bundle is simple, and
\begin{equation}\label{eq:specialExt}
\begin{aligned}
&\Ext^i_X(\cF,\cF)\simeq\CC\quad(i=0,1,2),&&\Ext^i_X(\cF,\cF)=0\quad(i\ge3),\\
&\Ext^1_X(\cF,\cF)\simeq H^1(\PP^2,\Omega^1_{\PP^2}),&&\Ext^2_X(\cF,\cF)\simeq H^2(\PP^2,\Omega^2_{\PP^2}).
\end{aligned}
\end{equation}

To relate spinor deformations to deformations of $X$, we use spinor modifications. Kuznetsov's modification theorem \cite[Theorem~1.2]{Kuz25} associates a conic bundle with an abstract spinor bundle while keeping the Clifford category fixed. At first order we can make this completely explicit. A first-order deformation of $\cF$ deforms its endomorphism algebra $\cB_0$ by commutators, and commutators with odd elements of the Clifford algebra are vector fields on the relative quadric; this produces a first-order deformation of $X$, together with a deformation of the Clifford category that is trivial up to Morita equivalence.

\begin{lemma}\label{lem:spinormod}
First-order spinor modification defines a natural linear map
\begin{equation}\label{eq:muf}
\mu_f\colon\Ext^1_X(\cF,\cF)\longrightarrow H^1(X,T_X),
\end{equation}
which is an isomorphism onto $H^1(X,T_X)^-\simeq H^1(Y,T_Y(-1,-1))$. The corresponding deformations of $X$ carry base-linear deformations of the Clifford category that are Morita trivial.
\end{lemma}

\begin{proof}
\emph{Construction.} Choose an affine cover of $\PP^2$. By \eqref{eq:Cliffordext}, a first-order deformation of $\cF$ is given by transition corrections $1+\epsilon a_{ij}$, where $(a_{ij})$ is a \v{C}ech cocycle with values in $\cB_0$, and the endomorphism algebra changes by $r\mapsto(1+\epsilon a_{ij})r(1-\epsilon a_{ij})=r+\epsilon[a_{ij},r]$. The scalar component acts trivially, and by \eqref{eq:specialExt} every class has a representative in the odd summand $\Omega^1_{\PP^2}$. Let $e_0$ be the generator of the summand $\cO_{\PP^2}$ of $\cE_0$, and use the conventions $uv+vu=2b(u,v)$ and $e_0^2=-1$, identifying $a\in\Omega^1_{\PP^2}$ with $ae_0$ in the even Clifford algebra. Commutation with $ae_0$ induces the vector field
\begin{equation}\label{eq:deltaa}
\delta_a(u)=-2za,\qquad\delta_a(z)=-2b(a,u),
\end{equation}
which preserves \eqref{eq:relquadric}, since $\delta_a\bigl(b(u,u)-z^2\bigr)=-4z\,b(u,a)+4z\,b(a,u)=0$. This polynomial identity holds over the dual numbers and also on degenerate conics. So the local automorphisms $1+\epsilon\delta_{a_{ij}}$ glue the local copies of $X$ to a first-order deformation of the conic bundle; the cocycle condition gives the descent, and changing the representative changes the result by an isomorphism. This defines $\mu_f$ without applying a theorem with a reduced-base hypothesis to $\PP^2\times\Spec\CC[\epsilon]/(\epsilon^2)$. The Clifford algebra of the deformed conic bundle is the deformed endomorphism algebra, which is Morita equivalent to the constant one through the deformed spinor module.

\emph{Identification.} The deformations obtained in this way are anti-invariant, and we compute them through the odd part of the tangent pushforward, $(\pi_*T_X)^-=T_Y\otimes\cO_Y(-1,-1)$. Twist the relative tangent sequence of $p$ by $\cO_Y(-1,-1)$. Its quotient $p^*T_{\PP^2}\otimes\cO_Y(-1,-1)$ has vanishing derived pushforward, being $\cO_{\PP^1}(-1)^{\oplus2}$ on every fibre, while $T_{Y/\PP^2}\otimes\cO_Y(-1,-1)\simeq\cO_Y(-2,1)$. Hence
\begin{equation}\label{eq:oddtangent}
Rp_*\bigl(T_Y\otimes\cO_Y(-1,-1)\bigr)\simeq p_*\cO_Y(-2,1)\simeq T_{\PP^2}(-3)\simeq\Omega^1_{\PP^2},
\end{equation}
all complexes being concentrated in degree zero. Under \eqref{eq:oddtangent}, the odd symbol of \eqref{eq:deltaa} is $-2a$, so the sheaf map underlying $\mu_f$ is an isomorphism from $\Omega^1_{\PP^2}$ onto the odd tangent pushforward. Taking $H^1(\PP^2,-)$ gives the claim.
\end{proof}

We now interpret the kernel of the differential of the period map in the special case in terms of the spinor bundle.

\begin{proposition}\label{prop:spinorkernel}
For a special Verra threefold, first-order spinor modification gives an isomorphism
\[
T_{[\cF]}\mathrm{Spin}_f(X)\xrightarrow{\ \sim\ }\Ker(\eta_f)=\Ker(d\cP)
\]
between one-dimensional spaces.
\end{proposition}

\begin{proof}
Lemma~\ref{lem:spinormod} identifies $T_{[\cF]}\mathrm{Spin}_f(X)$ with $H^1(X,T_X)^-$, which equals $\Ker(\eta_f)$ by Theorem~\ref{thm:Verracategorical} and $\Ker(d\cP)$ by Proposition~\ref{prop:specialVerra}. Equivalently, the modification gives a Morita-trivial deformation of the Clifford component, and the categorical computation identifies it with the kernel of the differential of the period map.
\end{proof}

The direction can also be described through the ambient bundle of the conic bundle. After twisting by a line bundle from the base, this ambient bundle is $\Omega^1_{\PP^2}(1)\oplus\cO_{\PP^2}(1)$ in the special case and $\cO_{\PP^2}^{\oplus3}$ in the ordinary case, and the spinor deformation class induces the Euler extension
\begin{equation}\label{eq:Euler}
0\longrightarrow\Omega^1_{\PP^2}(1)\longrightarrow\cO_{\PP^2}^{\oplus3}\longrightarrow\cO_{\PP^2}(1)\longrightarrow0,
\end{equation}
whose middle term is the ambient bundle of the ordinary case. Precisely, its class is a nonzero multiple of the generator of $\Ext^1(\cO_{\PP^2}(1),\Omega^1_{\PP^2}(1))$: the component $\Omega^1_{\PP^2}\to\Omega^2_{\PP^2}$ of the deformed Clifford transition map is a nonzero multiple of $v\mapsto a\wedge v$, and dualizing the trace-free Clifford bundle with the normalization of \cite[Remark~2.2]{Kuz25} gives the stated class. This describes the first-order deformation normal to the special locus; it does not mean that the fixed threefold carries a positive-dimensional family of spinor modifications.

\subsubsection{The obstruction}
In fact the first-order spinor deformation does not extend to second order.

\begin{proposition}\label{prop:obstruction}
For a special Verra threefold and $0\ne\alpha\in\Ext^1_X(\cF,\cF)$, the Yoneda square $\alpha\circ\alpha\in\Ext^2_X(\cF,\cF)$ is nonzero. Consequently the nonzero spinor deformation is obstructed at second order, and the completed miniversal deformation ring of $\cF$, with scalar automorphisms rigidified, is $\CC[[t]]/(t^2)$.
\end{proposition}

\begin{proof}
In the decomposition \eqref{eq:specialB0}, the component $\Omega^1_{\PP^2}\otimes\Omega^1_{\PP^2}\to\Omega^2_{\PP^2}$ of Clifford multiplication is a nonzero multiple of the wedge product. This follows from \eqref{eq:relquadric}: both odd Clifford factors contain the vector $e_0$, whose quadratic value is $-1$; see also \cite[Section~2.1, equation~(5) and Example~2.3]{Kuz25}. Under \eqref{eq:specialExt}, the Yoneda product is therefore a nonzero multiple of the cup product $H^1(\PP^2,\Omega^1_{\PP^2})^{\otimes2}\to H^2(\PP^2,\Omega^2_{\PP^2})$, and the square of the hyperplane class generates the target.

Since $\cF$ is a vector bundle whose only automorphisms are scalars, the Kuranishi construction of \cite[Lemma~2.4 and Section~2.3]{CPZ24} applies without any stability assumption, and its quadratic term is the Yoneda square, the primary obstruction to extending a first-order deformation to $\CC[t]/(t^3)$. So the miniversal deformation ring $\CC[[t]]/I$ has embedding dimension one, and $I\subset(t^2)$ contains an element of order two, which forces $I=(t^2)$.
\end{proof}

\subsection{Comparison of ordinary and special Verra threefolds}\label{ss:Verracompare}
Table~\ref{tab:Verra} collects the results of Sections~\ref{ss:Verraordinary} and~\ref{ss:Verraspecial}.

\begin{table}[t]
\centering
\renewcommand{\arraystretch}{1.3}
\begin{tabular}{>{\raggedright\arraybackslash}p{4.3cm}>{\raggedright\arraybackslash}p{4.4cm}>{\raggedright\arraybackslash}p{4.8cm}}
\toprule
 & ordinary Verra threefold & special Verra threefold\\
\midrule
model & divisor of bidegree $(2,2)$ in $\PP^2\times\PP^2$ & double cover of a $(1,1)$-divisor $Y$, branched along a K3 surface $S\in|-K_Y|$\\
$h^1(X,T_X)$ & $19$ & $19=18+1$ (invariant $+$ anti-invariant)\\
$\Ker(d\cP)$ & $0$ (Proposition~\ref{prop:VerraTorelli}) & $H^1(Y,T_Y(-1,-1))\simeq\CC$ (Proposition~\ref{prop:specialVerra})\\
$\Ker\mu$, $\NHH^2(\cE,X)$ & $0$, $0$ & $\CC r$, $\CC^2$\\
$H^0(X,\wedge^2T_X)$ & $0$ & $\CC$, anti-invariant\\
Clifford component $\cC_f$ & $\Db(\PP^2,\cO\oplus\cO(-2)^{\oplus3})$ & Enriques category, $S_{\cC_f}=\iota[2]$, $\HH^1(\cC_f)=0$\\
$\Ker(\eta_f)$ & $0$ & $H^1(X,T_X)^-\simeq\CC$\\
$\Ext^\bullet(\cF,\cF)$ & $\CC$ in degree $0$ & $\CC$ in degrees $0,1,2$\\
local moduli of $\cF$ & reduced point & $\Spec\CC[t]/(t^2)$\\
\bottomrule
\end{tabular}
\caption{Ordinary and special Verra threefolds, for a fixed conic bundle projection $f$.}
\label{tab:Verra}
\end{table}

For ordinary Verra threefolds, $h^1(X,T_X)=19$ is the dimension of the moduli space, $\dim|\cO_A(2,2)|-\dim\Aut(A)=35-16$; for special ones it was computed in Lemma~\ref{lem:Vdecomp}. We discuss the comparison under three headings.

\subsubsection{Moduli and the hyperelliptic analogy}
Special Verra threefolds form an $18$-dimensional family, a divisor in the $19$-dimensional moduli space of Verra threefolds. The covering involution forces the one-dimensional kernel $H^1(Y,T_Y(-1,-1))$, which is the direction normal to this divisor. This is closer to the situation of hyperelliptic curves than the Gushel--Mukai case, where the special locus has codimension three and the kernel grows from dimension two to dimension three.

\subsubsection{Kernels of the differentials of the classical and categorical period maps}
In both cases the kernel of the infinitesimal period map equals the kernel of the categorical period map of the Clifford component. For ordinary Verra threefolds this is immediate, since the exceptional collection has no normal Hochschild cohomology. For special ones, normal Hochschild cohomology is two-dimensional: besides the Kodaira--Spencer class of the anti-invariant deformation it contains the class of the bivector field, and the kernel of the differential of the period map is obtained by intersecting with the geometric summand of $\HH^2(X)$. The vanishing of $\HH^1(\cC_f)$, which in the special case comes from the Enriques structure of the Clifford component, plays the same role as the vanishing of $\HH^1(\Ku(X))$ for prime Fano threefolds of genus six and eight.

\subsubsection{Spinor bundles}
In both cases the kernel of the differential of the period map is the tangent space to the deformations of the canonical spinor bundle of a fixed conic bundle projection. The even Clifford algebra changes from $\cO\oplus\cO(-2)^{\oplus3}$, which has no higher cohomology, to $\cO\oplus\Omega^1\oplus\Omega^2$; its odd summand $\Omega^1_{\PP^2}$ carries the anti-invariant deformations, and its multiplication $\Omega^1\otimes\Omega^1\to\Omega^2$ is the wedge product. In the ordinary case the canonical spinor bundle is an isolated reduced point of its moduli space. In the special case it is an isolated nonreduced point, with a one-dimensional tangent space and a nonzero quadratic obstruction. So the kernel of the differential of the period map of a special Verra threefold is not the tangent space to a curve of pairwise distinct spinor bundles on the fixed threefold: the first-order spinor modification exists and deforms $X$ out of the special locus, but it does not extend to second order on $X$.


\begin{thebibliography}{BBF{\etalchar{+}}24}

\bibitem[AGG24]{AGG24}
R.~Aguilar, M.~Green, and P.~Griffiths, \emph{Infinitesimal invariants of mixed Hodge structures}, preprint, \texttt{arXiv:2406.17118}, 2024.

\bibitem[BBF{\etalchar{+}}24]{BBF+24}
A.~Bayer, S.~V.~Beentjes, S.~Feyzbakhsh, G.~Hein, D.~Martinelli, F.~Rezaee, and B.~Schmidt, \emph{The desingularization of the theta divisor of a cubic threefold as a moduli space}, Geom. Topol. \textbf{28} (2024), 127--160.

\bibitem[BF11]{BF11}
M.~C.~Brambilla and D.~Faenzi, \emph{Moduli spaces of rank-2 ACM bundles on prime Fano threefolds}, Michigan Math. J. \textbf{60} (2011), 113--148.

\bibitem[BFT23]{BFT23}
P.~Belmans, E.~Fatighenti, and F.~Tanturri, \emph{Polyvector fields for Fano 3-folds}, Math. Z. \textbf{304} (2023), Paper No.~12.

\bibitem[BP23]{BP23}
A.~Bayer and A.~Perry, \emph{Kuznetsov's Fano threefold conjecture via K3 categories and enhanced group actions}, J. Reine Angew. Math. \textbf{800} (2023), 107--153.

\bibitem[Cat84]{Cat84}
F.~Catanese, \emph{Infinitesimal Torelli theorems and counterexamples to Torelli problems}, in: Topics in Transcendental Algebraic Geometry (Princeton, 1981/1982), Ann. of Math. Stud. \textbf{106}, Princeton University Press, Princeton, NJ, 1984, 143--156.

\bibitem[CPZ24]{CPZ24}
H.~Chen, L.~Pertusi, and X.~Zhao, \emph{Some remarks about deformation theory and formality conjecture}, Ann. Univ. Ferrara \textbf{70} (2024), 761--779.

\bibitem[Deb20]{Deb20}
O.~Debarre, \emph{Gushel--Mukai varieties}, preprint, \texttt{arXiv:2001.03485}, 2020.

\bibitem[DIM12]{DIM12}
O.~Debarre, A.~Iliev, and L.~Manivel, \emph{On the period map for prime Fano threefolds of degree 10}, J. Algebraic Geom. \textbf{21} (2012), 21--59.

\bibitem[DK18]{DK18}
O.~Debarre and A.~Kuznetsov, \emph{Gushel--Mukai varieties: classification and birationalities}, Algebr. Geom. \textbf{5} (2018), 15--76.

\bibitem[Fae05]{Fae05}
D.~Faenzi, \emph{Bundles over the Fano threefold $V_5$}, Comm. Algebra \textbf{33} (2005), 3061--3080.

\bibitem[Fle86]{Fle86}
H.~Flenner, \emph{The infinitesimal Torelli problem for zero sets of sections of vector bundles}, Math. Z. \textbf{193} (1986), 307--322.

\bibitem[GP00]{GP00}
F.~J.~Gallego and B.~P.~Purnaprajna, \emph{Vanishing theorems and syzygies for K3 surfaces and Fano varieties}, J. Pure Appl. Algebra \textbf{146} (2000), 251--265.

\bibitem[Har10]{Har10}
R.~Hartshorne, \emph{Deformation Theory}, Graduate Texts in Mathematics \textbf{257}, Springer, New York, 2010.

\bibitem[HL10]{HL10}
D.~Huybrechts and M.~Lehn, \emph{The Geometry of Moduli Spaces of Sheaves}, 2nd ed., Cambridge Mathematical Library, Cambridge University Press, Cambridge, 2010.

\bibitem[Ili94]{Ili94}
A.~Iliev, \emph{The Fano surface of the Gushel threefold}, Compos. Math. \textbf{94} (1994), 81--107.

\bibitem[Ili97]{Ili97}
A.~Iliev, \emph{The theta divisor of bidegree $(2,2)$ threefold in $\PP^2\times\PP^2$}, Pacific J. Math. \textbf{180} (1997), 57--88.

\bibitem[JLLZ23]{JLLZ23}
A.~Jacovskis, X.~Lin, Z.~Liu, and S.~Zhang, \emph{Infinitesimal categorical Torelli theorems for Fano threefolds}, J. Pure Appl. Algebra \textbf{227} (2023), Paper No.~107418.

\bibitem[JLLZ24]{JLLZ24}
A.~Jacovskis, X.~Lin, Z.~Liu, and S.~Zhang, \emph{Categorical Torelli theorems for Gushel--Mukai threefolds}, J. Lond. Math. Soc. (2) \textbf{109} (2024), Paper No.~e12878.

\bibitem[JLZ22a]{JLZ22a}
A.~Jacovskis, Z.~Liu, and S.~Zhang, \emph{Brill--Noether theory for Kuznetsov components and refined categorical Torelli theorems for index one Fano threefolds}, preprint, \texttt{arXiv:2207.01021}, 2022.

\bibitem[Kon85]{Kon85}
K.~Konno, \emph{On deformations and the local Torelli problem of cyclic branched coverings}, Math. Ann. \textbf{271} (1985), 601--617.

\bibitem[KP17]{KP17}
A.~Kuznetsov and A.~Perry, \emph{Derived categories of cyclic covers and their branch divisors}, Selecta Math. (N.S.) \textbf{23} (2017), 389--423.

\bibitem[KP18]{KP18}
A.~Kuznetsov and A.~Perry, \emph{Derived categories of Gushel--Mukai varieties}, Compos. Math. \textbf{154} (2018), 1362--1406.

\bibitem[KP23]{KP23}
A.~Kuznetsov and A.~Perry, \emph{Categorical cones and quadratic homological projective duality}, Ann. Sci. \'Ec. Norm. Sup\'er. (4) \textbf{56} (2023), 1--57.

\bibitem[KPS18]{KPS18}
A.~G.~Kuznetsov, Yu.~G.~Prokhorov, and C.~A.~Shramov, \emph{Hilbert schemes of lines and conics and automorphism groups of Fano threefolds}, Jpn. J. Math. \textbf{13} (2018), 109--185.

\bibitem[Kuz04]{Kuz04}
A.~G.~Kuznetsov, \emph{Derived category of a cubic threefold and the variety $V_{14}$}, Tr. Mat. Inst. Steklova \textbf{246} (2004), 183--207.

\bibitem[Kuz06]{Kuz06}
A.~G.~Kuznetsov, \emph{Hyperplane sections and derived categories}, Izv. Math. \textbf{70} (2006), 447--547.

\bibitem[Kuz08]{Kuz08}
A.~Kuznetsov, \emph{Derived categories of quadric fibrations and intersections of quadrics}, Adv. Math. \textbf{218} (2008), 1340--1369.

\bibitem[Kuz09a]{Kuz09a}
A.~Kuznetsov, \emph{Hochschild homology and semiorthogonal decompositions}, preprint, \texttt{arXiv:0904.4330}, 2009.

\bibitem[Kuz09b]{Kuz09b}
A.~G.~Kuznetsov, \emph{Derived categories of Fano threefolds}, Tr. Mat. Inst. Steklova \textbf{264} (2009), 116--128.

\bibitem[Kuz15]{Kuz15}
A.~Kuznetsov, \emph{Height of exceptional collections and Hochschild cohomology of quasiphantom categories}, J. Reine Angew. Math. \textbf{708} (2015), 213--243.

\bibitem[Kuz25]{Kuz25}
A.~Kuznetsov, \emph{Spinor modifications of conic bundles and derived categories of $1$-nodal Fano threefolds}, preprint, \texttt{arXiv:2502.02082v2}, 2025.

\bibitem[Li19]{Li19}
C.~Li, \emph{Stability conditions on Fano threefolds of Picard number 1}, J. Eur. Math. Soc. (JEMS) \textbf{21} (2019), 709--726.

\bibitem[Log12]{Log12}
D.~Logachev, \emph{Fano threefolds of genus 6}, Asian J. Math. \textbf{16} (2012), 515--559.

\bibitem[LPZ23]{LPZ23}
C.~Li, L.~Pertusi, and X.~Zhao, \emph{Derived categories of hearts on Kuznetsov components}, J. Lond. Math. Soc. (2) \textbf{108} (2023), 2146--2174.

\bibitem[LZ22]{LZ22}
Z.~Liu and S.~Zhang, \emph{A note on Bridgeland moduli spaces and moduli spaces of sheaves on $X_{14}$ and $Y_3$}, Math. Z. \textbf{302} (2022), 803--837.

\bibitem[Per21]{Per21}
A.~Perry, \emph{Hochschild cohomology and group actions}, Math. Z. \textbf{297} (2021), 1273--1292.

\bibitem[PPZ26]{PPZ26}
A.~Perry, L.~Pertusi, and X.~Zhao, \emph{Moduli spaces of stable objects in Enriques categories}, Selecta Math. (N.S.) \textbf{32} (2026), Article~47.

\bibitem[Sno86]{Sno86}
D.~M.~Snow, \emph{Cohomology of twisted holomorphic forms on Grassmann manifolds and quadric hypersurfaces}, Math. Ann. \textbf{276} (1986), 159--176.

\bibitem[Usu81]{Usu81}
S.~Usui, \emph{Effect of automorphisms on variation of Hodge structures}, J. Math. Kyoto Univ. \textbf{21} (1981), 645--672.

\bibitem[Ver04]{Ver04}
A.~Verra, \emph{The Prym map has degree two on plane sextics}, in: The Fano Conference, Univ. Torino, Turin, 2004, 735--759.

\end{thebibliography}
\end{document}